\documentclass[mathpazo]{cicp}
\usepackage{lineno}
\usepackage{graphicx}
\usepackage{subfigure}
\modulolinenumbers[5]
\usepackage{float}

\usepackage{mathrsfs}
\usepackage{epstopdf}
\usepackage{amsmath}
\numberwithin{equation}{section}

\usepackage{amssymb}
\usepackage{latexsym}
\usepackage{color}
\usepackage{bm}
\usepackage{multicol}
\usepackage{multirow}
\usepackage{times}
\usepackage{amsthm}

\usepackage{amscd}
\usepackage{amsfonts}
\usepackage{ulem}

\usepackage{url}
\usepackage{dsfont}

\usepackage{mathrsfs}

\usepackage{booktabs}

\usepackage{hyperref}

\usepackage{verbatim}
\usepackage[latin1]{inputenc}
\usepackage{colortbl,dcolumn}
\usepackage{psfrag}

\newcommand{\E}{\mathbb E}
\newcommand{\R}{\mathbb R}

\newcommand{\AAA}{\mathcal A}

\graphicspath{ {fig/} {} }

\newtheorem{tm}{Theorem}[section]

\newtheorem{lm}{Lemma}[section]

\newtheorem{rk}{Remark}[section]

\allowdisplaybreaks[4]

\begin{document}

\title{Structure-preserving numerical methods for stochastic Poisson systems}


\author[Hong J L et.~al.]{Jialin Hong\affil{1},
        Jialin Ruan\affil{2}, Liying Sun\affil{1}, Lijin Wang\affil{2}\comma\corrauth}
 \address{\affilnum{2}\ School of Mathematical Sciences,
           University of Chinese Academy of Sciences, 19 YuQuan Road, Shijingshan District, 
          Beijing 100049, China\\
          \affilnum{1}\ LSEC, ICMSEC,
Academy of Mathematics and Systems Science, Chinese Academy of Sciences, Beijing,  100190, China \\
             }
 \emails{{\tt hjl@lsec.cc.ac.cn} (J.~Hong), {\tt rjl2011@mail.ustc.edu.cn} (J.~Ruan), {\tt liyingsun@lsec.cc.ac.cn} (L.~Sun), {\tt ljwang@ucas.ac.cn} (L.~Wang)}

\begin{abstract}
We propose a class of numerical integration methods for stochastic Poisson systems (SPSs) of arbitrary dimensions. Based on the Darboux-Lie theorem, we transform the  SPSs to their canonical form, the generalized stochastic Hamiltonian systems (SHSs), via canonical coordinate transformations found by solving certain PDEs defined by the Poisson brackets of the SPSs. An $\alpha$-generating function approach with $\alpha\in[0,1]$ is then used to create symplectic discretizations of the SHSs, which are then transformed back by the inverse coordinate transformation to numerical integrators for the SPSs. These integrators are proved to preserve both the Poisson structure and the Casimir functions of the SPSs. Applications to a three-dimensional stochastic rigid body system and a three-dimensional stochastic Lotka-Volterra system show efficiency of the proposed methods.
\end{abstract}

\ams{60H35, 60H15, 65C30, 60H10, 65D30}
\keywords{stochastic Poisson systems, Poisson structure, Casimir functions, Poisson integrators, symplectic integrators, generating functions, stochastic rigid body system, stochastic Lotka-Volterra system}

\maketitle

\section{Introduction}
Poisson systems form a class of important mechanical systems whose long history dates back to the 19th century (\cite{Darboux,Lie,ECG2002}).  As a generalizatioin of the Hamiltonian systems which are defined on even-dimensional symplectic manifolds, the poisson systems possess similar but extended structural properties, and can be defined on Poisson manifolds of arbitrary dimensions. They have a large scope of applications, such as in astronomy, robotics, fluid mechanics, electrodynamics, quantum mechanics, nonlinear waves, and so on (\cite{ZQ1994}). Unlike Hamiltonian systems where plenty literatures are available on their numerical approximations, there have not been as many studies on numerical simulations for the Poisson systems. One of the main challenges for numerical approximations of the Poisson systems is that such approximations depend on the concrete structure matrix, which makes it difficult to establish general methodologies (\cite{fengqbook,ECG2002}).

Symplectic methods for Hamiltonian systems have been developed during the last decades (\cite{fengqbook,ECG2002,qwbook}). They find applications in many fields where Hamiltonian systems appear, and are proved to be much superior than non-symplectic methods in long time simulation, due to their ability of preserving the symplectic structure of the original systems (see e.g. \cite{mr2,mr3,mr7}). Structure-preserving algorithms of a broader sense are then aroused which seek for preservation of more structural conservation law in numerical discretization, such as energy, momentum, etc. (see e.g. \cite{mr1,mr4,mr6}). The Poisson structure is an extension of the symplectic structure to arbitrary-dimension and variable structure matrices, and is reduced to the symplectic structure when the structure matrices degenerate to the even-dimensional symplectic matrix $J$. It is an intrinsic structure of the Poisson systems. However, it has been observed that, symplectic methods in general do not preserve the Poisson structure (\cite{fengqbook,ECG2002,SQ2005}). Therefore, there is a need to develop Poisson integrators which can inherit the Poisson structure of the Poisson systems. Such attempts have been made for deterministic cases in e.g. \cite{cohen1,ge,62,63,65,71,73,ZQ1994} etc. 

In recent years, there arise some numerical studies on certain special stochastic Poisson systems (SPSs). \cite{cohen} proposed a class of energy-preserving numerical methods for stochastic Poisson systems where the deterministic and stochastic Hamiltonians vary by a constant. These methods are proved to preserve quadratic Casimir functions as well. \cite{mading} constructed a class of explicit parametric stochastic Runge-Kutta methods with truncated random variables for such stochastic Poisson systems, and showed that these methods are energy-preserving for suitable parameters, and can be of any prescribed convergence orders. For stochastic Poisson systems of even dimensions and invertible structure matrices, \cite{xu} investigates structure-preserving Runge-Kutta and partitioned Runge-Kutta type methods. Up to now, we have not seen numerical analysis on structure-preserving algorithms for general stochastic Poisson systems with arbitrary dimensions, different Hamiltonians and multiple noises. 

In this paper, we propose a class of numerical methods for general stochastic Poisson systems. By appropriate coordinate transformations, we rewrite the SPSs into their canonical forms, which are generalized stochastic Hamiltonian systems (SHSs). Then we apply a stochastic $\alpha$-generating function approach to construct symplectic schemes for the resulted SHSs, and transform the symplectic schemes back to numerical schemes for the SPSs afterwards.  The so-proposed methods are shown to preserve the Poisson structure and the Casimir functions of the SPSs. Suitable coordinate transformations are found by solving certain partial differential equations.  As applications, we apply the proposed methods to a three-dimensional stochastic rigid body system and a three-dimensional stochastic Lotka-Volterra system. 

Contents of the paper are organized as follows. In Section 2 we introduce the concept of stochastic Poisson systems, and prove that they possess the Poisson structure, and the Casimir functions are invariant quantities of the SPSs.  In Section 3 we use the $\alpha$-generating function approach and the Darboux-Lie theorem to construct numerical methods for the SPSs, and apply them to the stochastic rigid body system and the stochastic Lotka-Volterra system.  Numerical experiments are illustrated in section 4, followed by a few concluding remarks in Section 5.

\section{The stochastic Poisson systems (SPSs)}
\label{sec;SPS}

Consider the following $d$-dimensional stochastic Poisson system 
\begin{align}
\label{sps-eq-1}
d y &= B(y) \left(\nabla K_0(y) d t +
\sum_{r=1}^m\nabla K_r(y) \circ d W_r(t)\right),\nonumber\\
y(t_0)&=y_0,
\end{align}
where $t\in[t_0,T]$, $y\in \R^d$,  $W(t)=(W_1(t),\cdots,W_m(t))$ is an m-dimensional standard Wiener process defined on a complete filtered probability space $(\Omega, \mathcal{F}, \mathcal{P}, \{\mathcal{F}_t\}_{t\ge 0})$, and the symbol $``\circ"$ represents the Stratonovich product. $y_0$ satisfies
\begin{align}\label{y0}
&(a)\,\,\, \,E|y_0|^2<\infty, \,\,\,\mbox{with} |\cdot| \,\,\,\mbox{being the Euclidean norm, and} \nonumber\\
&(b)\,\,\, \,
y_0 \,\,\,\,\,\mbox{is}\,\,\,\mathcal{F}_0-\mbox{measurable}.
\end{align}
$B:\R^d\mapsto\R^{d\times d}$  and $K_i: \R^d\mapsto\R$ $(i=0,\cdots,m)$  are sufficiently smooth functions, and we assume the coefficients 
$$a(y):=B(y)\nabla K_0(y),\qquad  b_r(y):=B(y)\nabla K_r(y)\,\,\,(r=1,\ldots,m)$$
satisfy the conditions guaranteeing existence and uniqueness of the solution of the stochastic differential equations system (\ref{sps-eq-1}) (see e.g. \cite{oksendal}), namely,
\begin{align}\label{addlip}
|a(x)-a(y)|+\sum_{r=1}^m|b_r(x)-b_r(y)|\le D_1|x-y|,\quad x,\,\,y\in \R^d,\,\,\,\mbox{for some }\,\,D_1>0;\\
|a(y)|+\sum_{r=1}^m |b_r(y)|\le D_2(1+|y|),\quad y\in \R^d,\,\,\,\mbox{for some}\,\,D_2>0.\label{addlip2}
\end{align}
It is also known that (see e.g. \cite{oksendal}), under the conditions (\ref{y0})-(\ref{addlip2}), the solution $y(t,y_0,\omega)$ of (\ref{sps-eq-1}) will not blow up in finite time interval $[t_0,T]$, namely $E[\int_{t_0}^T|y(t)|^2dt]<\infty$.

Further, $B(y)=\big( b_{i j}(y) \big)\in \R^{d\times d}$ is skew-symmetric, that is,
\begin{equation}\label{as2}
b_{i j}(y) =-b_{j i}(y),
\end{equation} and satisfies the condition
\begin{equation}
\label{dps-eq-1}
\sum\limits_{s=1}^{d} \left( \frac{\partial b_{i j}(y)}{\partial y_s} b_{s k}(y) + \frac{\partial b_{j k}(y)}{\partial y_s} b_{s i}(y) + \frac{\partial b_{k i}(y)}{\partial y_s} b_{s j}(y) \right) = 0, 
\end{equation}
for all $i,j,k$. In addition, we assume that $B(y)$ is of constant rank $d-l=2n$ with $l\ge 0$.
The SPS \eqref{sps-eq-1} will degenerate to the stochastic Hamiltonian system (SHS) (\cite{mil1,mil2}) when the matrix $B(y)=J^{-1}=\begin{bmatrix}
O& -I_n\\
I_n & O
\end{bmatrix}$ with $d=2n$. 
\subsection{The Poisson structure and Casimir functions of the SPSs}
As was given in \cite{ECG2002} for deterministic cases, the structure matrix $B(y)$ characterized by the properties (\ref{as2} )and (\ref{dps-eq-1}) defines the Poisson bracket $\{F,G\}$ of two smooth functions $F(y)$ and $G(y)$ as
\begin{equation}
\{ F, G \}(y) := \sum\limits_{i, j=1}^{n}
\frac{\partial F(y)}{\partial y_i} b_{i j}(y) \frac{\partial G(y)}{\partial y_j},
\end{equation}
or in vector notation 
\begin{equation}\label{vector}
\{ F, G \}(y) =  \nabla F(y)^{\mathrm{T}} B(y) \nabla G(y).
\end{equation}
The Poisson bracket$\{\cdot,\cdot\}$ is bilinear, skew-symmetric, 
and satisfies the Jacobi identity
\begin{equation*}
\{ \{ F, G \}, H \} + \{ \{ G, H \}, F \} + \{ \{ H, F \}, G \} = 0,
\end{equation*}
and the Leibniz rule
\begin{equation*}
\{ F \cdot G, H \} = F \cdot \{ G, H \} + G \cdot \{ H, F \}.
\end{equation*}
A map $\varphi : U \to \R^d$ (where $U$ is an open set in $\R^d$) is called a Poisson map if it commutes with the Poisson bracket, namely,
\begin{equation}\label{poissonmap}
	\{F\circ \varphi,G\circ\varphi\}(y)=\{F,G\}(\varphi(y)),
\end{equation}
for all smooth functions $F,G$ defined on $\varphi(U)$.
An identical expression of (\ref{poissonmap}) that we use in the following discussion is
\begin{equation}
\label{sps-eq-2}
\left[\frac{\partial \varphi(y)}{\partial y}\right]B(y)
\left[\frac{\partial \varphi(y)}{\partial y}\right]^\top
=B(\varphi(y)).
\end{equation}
The equivalence of (\ref{poissonmap}) and (\ref{sps-eq-2}) can be proved easily by using the vector formulation of the Poisson bracket (\ref{vector}) and the differential chain rule. 

\begin{lm} \label{lm-1}(\cite{Darboux,Lie,ECG2002}) Suppose that the matrix $B(y)$ defines a Poisson bracket and is of constant rank $d-l=2n$ in a neighborhood of $y_0\in\R^d$. Then there exist functions $P_1(y),\ldots,P_n(y),Q_1(y),\ldots,Q_n(y)$, and $C_1(y),\ldots,C_l(y)$ satisfying 
\begin{equation}\label{dltmeq} 
\begin{matrix}
&\{P_i,P_j\}=0, &\{P_i,Q_j\}=-\delta_{ij}, &\{P_i,C_s\}=0, \\
&\{Q_i,P_j\}=\delta_{ij}, &\{Q_i,Q_j\}=0, &\{Q_i,C_s\}=0,\\
&\{C_k,P_j\}=0, &\{C_k,Q_j\}=0, &\{C_k,C_s\}=0
\end{matrix}
\end{equation}
for $i=1,\ldots,n,j=1,\ldots,n,k=1,\ldots,l,s=1,\ldots,l$, on a neighborhood of $y_0$. The gradients of $P_i,\,\,Q_j,\,\,C_k$ $(i,=1,\ldots,n,\,\,j=1,\ldots,n,\,\,k=1,\ldots,l)$ are linearly independent, so that the $\R^d\rightarrow\R^d$ mapping $y\rightarrow(P_1(y),\ldots,P_n(y),Q_1(y),\ldots,Q_n(y),C_1(y),\ldots,C_l(y))$ constitutes a local change of coordinates to canonical form.
\end{lm}

Lemma \ref{lm-1} is also called the Darboux-Lie theorem. Next we use this theorem to prove the Poisson structure of the stochastic Poisson systems. 

\begin{tm}\label{tm-1} Under the conditions (\ref{y0})-(\ref{dps-eq-1}), for each $t$, almost surely, the solution flow $\varphi_t$ of the stochastic Poisson system (\ref{sps-eq-1})  is a Poisson map wherever it is defined.
\end{tm}
{\bf Proof.} 
Under the conditions (\ref{y0})-(\ref{addlip2}), almost surely, there exists an unique solution $y=y(y_0,t,\omega)$ of (\ref{sps-eq-1}), where $\omega\in \Omega_0\subset \Omega$ with $\mathcal{P}(\Omega_0)=1$.Due to (\ref{as2}) and (\ref{dps-eq-1}),  $B(y)$ of (\ref{sps-eq-1}) can define a Poisson bracket, and is of constant rank $d-l=2n$ by assumption. Then, according to the Darboux-Lie theorem, there exist functions $P_1(y),\ldots,P_n(y),Q_1(y),\ldots,Q_n(y)$, and $C_1(y),\ldots,C_l(y)$ satisfying (\ref{dltmeq}) such that the gradients of $P_i,\,\,Q_j,\,\,C_k$ $(i,=1,\ldots,n,\,\,j=1,\ldots,n,\,\,k=1,\ldots,l)$ are linearly independent, and the mapping $$y\rightarrow \bar{y}=:(P_1(y),\ldots,P_n(y),Q_1(y),\ldots,Q_n(y),C_1(y),\ldots,C_l(y)):=\theta(y)$$ constitutes a change of coordinates. The invertible Jacobian matrix of the coordinates transformation is 
$$ \frac{\partial \bar{y}}{\partial y}=\begin{pmatrix}\nabla P_1(y),\cdots,\nabla P_n(y),\nabla Q_1(y),\cdots \nabla Q_n(y),\nabla C_1(y),\cdots,\nabla C_l(y)\end{pmatrix}^\top=:A(y).$$ 
Then according to the vector representation of the Poisson bracket (\ref{vector}), as well as (\ref{dltmeq}), we have
\begin{align}\label{apple}
& \frac{\partial \bar{y}}{\partial y} B(y)\frac{\partial \bar{y}}{\partial y}^\top=A(y)B(y)A(y)^\top\nonumber\\&=\left(\begin{matrix}\{P_1,P_1\}&\cdots&\{P_1,P_n\}&\{P_1,Q_1\}&\cdots&\{P_1,Q_n\}&\{P_1,C_1\}&\cdots&\{P_1,C_l\}\\ \vdots&\cdots&\vdots&\vdots&\cdots&\vdots&\vdots&\cdots&\vdots\\\{P_n,P_1\}&\cdots&\{P_n,P_n\}&\{P_n,Q_1\}&\cdots&\{P_n,Q_n\}&\{P_n,C_1\}&\cdots&\{P_n,C_l\}\\
\{Q_1,P_1\}&\cdots&\{Q_1,P_n\}&\{Q_1,Q_1\}&\cdots&\{Q_1,Q_n\}&\{Q_1,C_1\}&\cdots&\{Q_1,C_l\}\\ \vdots&\cdots&\vdots&\vdots&\cdots&\vdots&\vdots&\cdots&\vdots\\\{Q_n,P_1\}&\cdots&\{Q_n,P_n\}&\{Q_n,Q_1\}&\cdots&\{Q_n,Q_n\}&\{Q_n,C_1\}&\cdots&\{Q_n,C_l\}\\ 
 \{C_1,P_1\}&\cdots&\{C_1,P_n\}&\{C_1,Q_1\}&\cdots&\{C_1,Q_n\}&\{C_1,C_1\}&\cdots&\{C_1,C_l\}\\ \vdots&\cdots&\vdots&\vdots&\cdots&\vdots&\vdots&\cdots&\vdots\\\{C_l,P_1\}&\cdots&\{C_l,P_n\}&\{C_l,Q_1\}&\cdots&\{C_l,Q_n\}&\{C_l,C_1\}&\cdots&\{C_l,C_l\} \end{matrix}\right)   \nonumber \\
 &=\left(\begin{matrix}0&-I&0\\ I&0&0\\ 0&0&0\end{matrix}\right)=\left(\begin{matrix}J^{-1}&0\\0&0\end{matrix}\right).
 \end{align}
 Then, for each $\omega\in \Omega_0$, it holds for $y=y(y_0,t,\omega)$ that
\begin{align}\label{consps}
d\bar{y}&=\frac{\partial \bar{y}}{\partial y}dy=\frac{\partial \bar{y}}{\partial y}\left(B(y)\left(\nabla K_0(y)dt+\sum_{r=1}^m\nabla K_r(y)\circ d W_r(t)\right)\right)\nonumber\\
&=\frac{\partial \bar{y}}{\partial y}\left(B(y)\left(\frac{\partial \bar{y}}{\partial y}^\top\nabla H_0(\bar{y})dt+\frac{\partial \bar{y}}{\partial y}^\top\sum_{r=1}^m\nabla H_r(\bar{y})\circ dW_r(t)\right)\right)\nonumber\\
&=\left(\begin{matrix}J^{-1}&0\\0&0\end{matrix}\right)\left(\nabla H_0(\bar{y})dt+\sum_{r=1}^m\nabla H_r(\bar{y})\circ d W_r(t)\right),
\end{align}
where 
 $H_i(\bar{y})=K_i(y)$ for $i=0,1\cdots,m$. Note that the Stratonovich chain rule is necessary for the validity of (\ref{consps}). The number of the zero rows in the structure matrix $\left(\begin{matrix}J^{-1}&0\\0&0\end{matrix}\right)$ is $l$, depending on the rank $2n$ of the matrix $B(y)$, since $2n+l=d$ by assumption. Denote $\bar{y}=(Z(y)^\top,C(y)^\top)^\top$, where $$Z(y)=(P(y)^\top, Q(y)^\top)^\top,\quad P(y)=(P_1(y),\cdots,P_n(y))^\top,$$ $$Q(y)=(Q_1(y),\cdots,Q_n(y))^\top, \quad C(y)=(C_1(y),\cdots,C_l(y))^\top,$$
then (\ref{consps}) is equivalent to 
\begin{eqnarray}\label{consps1}
dZ&=&J^{-1}\left(\nabla_Z H_0(Z,C)dt+\sum_{r=1}^m\nabla_Z H_r(Z,C)\circ dW_r(t)\right),\label{consps1a}\\
dC&=&0,\label{consps1b}
\end{eqnarray}
where the first equation is a $2n$-dimensional stochastic Hamiltonian system with constant parameters $C$. Due to the symplecticity of the SHSs (\cite{mil1,mil2,SW2017}), and with the help of the formulation (\ref{consps1a}), we can derive that, for each $t$, almost surely, the flow $\psi_t(\bar{y}_0)$ of the system (\ref{consps}) satisfies
\begin{equation}\label{gensym}
\frac{\partial \psi_t(\bar{y}_0)}{\partial \bar{y}_0}\left(\begin{matrix}J^{-1}&0\\0&0\end{matrix}\right)\frac{\partial \psi_t(\bar{y}_0)}{\partial \bar{y}_0}^\top=\left(\begin{matrix}J^{-1}&0\\0&0\end{matrix}\right),
\end{equation}
where $\bar{y}_0=\theta(y_0)$.
Moreover, since it holds $$\psi_t(\bar{y}_0)=\theta(\varphi_t(y_0)),$$ taking partial derivative with respect to $y_0$ on both sides of the above equation, we have
\begin{equation}\label{pf1}
\frac{\partial \psi_t(\bar{y}_0)}{\partial \bar{y}_0}=\theta'(\varphi_t(y_0))\frac{\partial \varphi_t(y_0)}{\partial y_0}A(y_0)^{-1}.
\end{equation}
Substituting (\ref{pf1}) into (\ref{gensym}), we get
\begin{equation}\label{pf2}
\theta'(\varphi_t(y_0))\frac{\partial \varphi_t(y_0)}{\partial y_0}B(y_0)\frac{\partial \varphi_t(y_0)}{\partial y_0}^\top\theta'(\varphi_t(y_0))^\top=\left(\begin{matrix}J^{-1}&0\\0&0\end{matrix}\right),
\end{equation}
due to (\ref{apple}). Note that $\theta'(\varphi_t(y_0))=A(\varphi_t(y_0))$, then with replacing $y$ in (\ref{apple}) by $\varphi_t(y_0)$, we obtain the equivalent form of (\ref{pf2})
\begin{equation}\label{pf3}
\frac{\partial \varphi_t(y_0)}{\partial y_0}B(y_0)\frac{\partial \varphi_t(y_0)}{\partial y_0}^\top=B(\varphi_t(y_0)), 
\end{equation}
almost surely, for any $t$ and any $y_0$, as long as the solution remains in the definition domain of $B$ and $K_i$ $(i=0,\dots,m)$.
This ends the proof. \hfill $\square$
\begin{remark}\label{rmadd}
We note that, the existence of $\frac{\partial \varphi_t(y_0)}{\partial y_0}$ can be assured by continuous differentiability of $B(\cdot)$ and twice continuous differentiability of $K_i(\cdot)$ $(i=0,\dots,m)$. Denote $z=\frac{\partial \varphi_t(y_0)}{\partial y_0}=(z_1,\dots,z_d)$, where $z_j=\frac{\partial \varphi_t(y_0)}{\partial y_0^j}$ ($j=1,\dots,d$) is the $j-th$ column of $z$. The SDE for $z_j$ obtained by formally taking partial derivative with respect to $y_0^j$ on both sides of (\ref{sps-eq-1}) is
\begin{align}\label{augm}
dz_j=M_0(t,y_0)z_jdt+\sum_{r=1}^m M_r(t,y_0)z_j\circ dW_r(t), \quad z_j(t_0)=e_j,
\end{align}
where $e_j$ is the $j$-th column of the $d$-dimensional identity matrix, $M_i(t,y_0)\in\R^{d\times d}$ ($i=0,\dots,m$) and
\begin{align}\label{zjeq}
M_i(t,y_0)&=B'(\varphi_t(y_0))(\nabla K_i(\varphi_t(y_0)))+B(\varphi_t(y_0))\nabla^2K_i(\varphi_t(y_0)).
\end{align}
From (\ref{zjeq}) we see that, for any given $y_0$, if $B(\cdot)$ is continuously differentiable and $K_i(\cdot)$ ($i=0,\dots,m$) are twice continuously differentiable, $M_i(t,y_0)$ are $\R^{d\times d}$-valued continuous functions of $t$, which are then bounded in norm on $[t_0,T]$. Thus the linear SDE (\ref{augm}) has a unique solution $z_j$ almost surely ($j=1,\dots,d$).

\end{remark}

We call (\ref{consps}) the canonical form of the SPS (\ref{sps-eq-1}), which is a generalized stochastic Hamiltonian system, and the transformation $\theta$ the canonical transformation.

Similar to the deterministic case (see e.g. \cite{ECG2002}), we call a function $C(y)$ a Casimir function of the system \eqref{sps-eq-1} if
\begin{equation}
\label{sps-eq-3}
\nabla C(y)^{\top} B(y) = 0 \quad \mbox{for all} \,\,\, y, \,\,\,\mbox{almost surely.}
\end{equation}
Thus, by (\ref{sps-eq-1}) and the Stratonovich chain rule, we have
\begin{align*}
dC(y)=\nabla C(y)^{\top} B(y)\left(\nabla K_0 (y)dt+\sum_{r=1}^m\nabla K_r(y)\circ dW_r(t)\right)=0.
\end{align*}
Therefore, each Casimir function $C(y)$ is a first integral of its corresponding SPS.  The existence and concrete forms of the Casimir functions depend merely on the structure matrix $B(y)$ of the SPS, whatever the Hamiltonians $K_i(y)$ $(i=0,\cdots,m)$ are. 

\subsection{SPSs in applications}\label{sec2.3}
Here we present two concrete models of stochastic Poisson systems in applications. 
\subsubsection{The stochastic rigid body system (SRB) (\cite{cohen,liao})} 
Consider the system
\begin{equation}\label{exa1}
dy=B_1 (y)\nabla K_1 (y)(dt+c_1\circ dW(t)),
\end{equation}
where $y=(y_1,y_2,y_3)^\top$, 
$K_1 (y) := \frac 12
\left(
\frac{y_1^2}{I_1}
+\frac{y_2^2}{I_2}
+\frac{y_3^2}{I_3}
\right)$,  $I_1$, $I_2$, $I_3$ and $c_1$ are constants and
$$B_1 (y) =
\begin{bmatrix}
0& -y_3 & y_2\\
y_3 & 0& -y_1\\
-y_2 & y_1& 0
\end{bmatrix}.$$
Since the matrix $B_1(y)$ satisfies the conditions \eqref{as2} and \eqref{dps-eq-1}, $B_1(\cdot) $ is continuously differentiable, $K_1(\cdot) $is twice continuously differentiable, \eqref{exa1} possesses the Poisson structure \eqref{pf3}.
The Casimir function of \eqref{exa1} is a quadratic function $$C_1(y):=\frac{1}{2}\left(y_1^2+y_2^2+y_3^2\right),$$ since $\nabla C_1 (y)^{\top} B_1 (y) = 0.$ Given the initial value $(y^0_1,y^0_2,y^0_3)^\top,$ the Casimir function $C_1(y)$ is always equal to the constant $\mathcal C_1=\frac{1}{2}\left((y_1^0)^2+(y_2^0)^2+(y_3^0)^2\right).$

The stochastic rigid body system has many physical applications. For instance, it can be used to describe the roll motion of the ship under impact of severe weather conditions (\cite{ACO2005})


\subsubsection{The stochastic Lotka-Volterra system (SLV) (\cite{cohen})} \label{lvexa}
Consider the system
\begin{equation}\label{exa2}
dy = B_2 (y) \nabla K_2 (y)(dt + c_2 \circ dW(t)),
\end{equation}
where $y=(y_1,y_2,y_3)^\top$, 
$K_2 (y) = ab y_1+y_2-a y_3+\nu \ln y_2 - \mu \ln y_3$
is a continuous function,
$$B_2 (y) =
\begin{bmatrix}
0& ry_1y_2& bry_1y_3\\
-ry_1y_2& 0& y_2y_3\\
-bry_1y_3& -y_2y_3& 0
\end{bmatrix},
$$ and $a$, $b$, $c_2$, $r$, $\nu$, $\mu$ are constants. Note that the system (\ref{exa2}) can be regarded as belonging to a generalization of the stochastic Lotka-Volterra systems (9) in \cite{MMR2001}, but of Stratonovich sense.

We can check the condition \eqref{sps-eq-3} to know that $C_2 (y) := \frac{1}{r}\ln y_1 - b \ln y_2 + \ln y_3 $ is the Casimir function of the stochastic Lotka-Volterra system.
Given the initial value $(y^0_1,y^0_2,y^0_3)^\top,$ the Casimir function $C_2(y)$ is always equal to the constant $\mathcal C_2= \frac{1}{r}\ln y_1^0 - b \ln y_2^0 + \ln y_3^0.$
Moreover, it can be verified that, with positive initial value $y^0$, the solution to \eqref{exa2} remains positive almost surely. 

The stochastic Lotka--Volterra system characterizes the population systems under stochastic influences which are inevitable and unignorable. For instance,  \cite{MMR2001} reveals that the environmental noises can suppress the potential population explosion.


\section{Numerical solution of the SPSs} 
\subsection{The transformation to canonical form}
By the Darboux-Lie theorem, a SPS  can be transformed to a generalized SHS. 
Our strategy is to first construct symplectic schemes for the generalized SHS, and then transform them back to obtain Poisson schemes for the SPS. The following are the details of this procedure:
\begin{itemize}
	\item  Given the $d$-dimension SPS \eqref{sps-eq-1}, find the Casimir functions $C_1(y) \equiv  \mathcal{C}_1, \dots $, $C_{l}(y)\equiv \mathcal{C}_l$  by solving $\nabla C(y)^\top B(y)=0$, where $l=d-2n$, and $2n$ is the rank of the anti-symmetric matrix $B(y)$. Denote $\mathcal C(y)=(C_1(y),\cdots,C_l(y))^\top$.
	\item 
	Use the  coordinate transformation $\bar{y}=\theta(y)$ as described in the proof of Theorem \ref{tm-1},  which has an invertible Jacobian matrix $A(y) =\theta'(y)= \frac{\partial \bar{y}}{\partial y},$ to transform the structure matrix $B(y)$ to $A(y) B(y) A(y)^{\mathrm{T}}.$
	
	\item Let  $A(y) B(y) A(y)^{\mathrm{T}} = B_0$ with $B_0$ being the constant structure matrix of a generalized SHS, e.g., $B_0 = \begin{bmatrix} J^{-1}_{2n\times2n} & 0_{2n\times(d-2n)} \\ 0_{(d-2n)\times2n} & 0_{(d-2n)\times(d-2n)} \end{bmatrix}$, to solve for the coordinate transformation 
\begin{equation}\label{ctc}
	\bar{y}=\theta(y)=(P_1(y),\cdots,P_n(y),Q_1(y),\cdots,Q_n(y),C_1(y),\cdots,C_l(y))^\top,  
\end{equation} 
and its inverse $y=\theta^{-1}(\bar{y})$. Note that the last $l$ coordinates in \eqref{ctc} are just the $l$ Casimir functions (by concrete calculations, or referring to the proof of the Darboux-Lie Theorem in e.g. \cite{ECG2002}).
Then we obtain the generalized SHS
\begin{equation}\label{GSHS}
d \bar{y} = B_0\left (\nabla H_0(\bar{y}) d t + \sum_{r=1}^m\nabla H_r(\bar{y}) \circ d W_r(t)\right ),
\end{equation}
where $H_i(\bar y)=K_i(y)$ $(i=0,1,\cdots,m)$. As discussed for (\ref{consps1a})-(\ref{consps1b}), if we denote $\bar{y}=(Z(y)^\top,C(y)^\top)^\top$, where $$Z(y)=(P(y)^\top, Q(y)^\top)^\top,\quad P(y)=(P_1(y),\cdots,P_n(y))^\top,$$ $$Q(y)=(Q_1(y),\cdots,Q_n(y))^\top, \quad C(y)=(C_1(y),\cdots,C_l(y))^\top,$$
then (\ref{GSHS}) is equivalent to 
\begin{eqnarray}\label{consps1}
dZ&=&J^{-1}\left(\nabla_Z H_0(Z,C)dt+\sum_{r=1}^m\nabla_Z H_r(Z,C)\circ dW_r(t)\right),\label{consps1a-new}\\
dC&=&0,\label{consps1b-new}
\end{eqnarray}
where (\ref{consps1a-new}) is a stochastic Hamiltonian system (SHS) with constant parameter vector $C$.

\item Given initial value $y_0$ of the SPS (\ref{sps-eq-1}), we can get $Z_0=(P(y_0)^\top,Q(y_0)^\top)^\top$. Apply a symplectic scheme $Z^{j+1}=\psi_h(Z^j,\mathcal{C})$ to the SHS \eqref{consps1a-new} with constant parameters $C$, where $Z^{j+1}=(P_1^{j+1},\ldots,P_n^{j+1},Q_1^{j+1},\ldots,Q_n^{j+1})^\top$ ($j=0,1,2,\cdots$) denotes the $(j+1)$-th step numerical value that approximates $Z(t_0+(j+1)h)$ of the solution of (\ref{consps1a-new}).  Then use the inverse transformation $y=\theta^{-1}(\bar{y})$ to transform $\bar{y}^{j+1}=((Z^{j+1})^\top,C^\top)^\top$ back to $y^{j+1}$, namely, $y^{j+1}=\theta^{-1}(\bar{y}^{j+1})$, to obtain the numerical value $y^{j+1}$ ($j=0,1,2,\cdots$) that approximates $y(t_0+(j+1)h)$ of the solution of (\ref{sps-eq-1}), which we denote by $y^{j+1}=\varphi_h(y^j)=\theta^{-1}(\bar{y}^{j+1})$, where $y=(y_1,\ldots,y_d)^\top$.
\end{itemize}
\begin{tm}\label{gen1} The above obtained numerical schemes $y^{j+1}=\varphi_h(y^j)$ are stochastic Poisson integrators for the stochastic Poisson system \eqref{sps-eq-1}, namely, they preserve both the Poisson structure and the Casimir functions of the SPS \eqref{sps-eq-1} almost surely.
\end{tm}
{\bf Proof.} 
Since $\bar y=\theta (y)=(P(y)^\top,Q(y)^\top, \mathcal C(y)^\top)^\top$, then
\begin{align}
 \psi_h(P^j,Q^j,\mathcal C)&=\theta(\varphi_h(y^j))=(P(\varphi_h(y^j))^\top,Q(\varphi_h(y^j))^\top, \mathcal C(\varphi_h(y^j))^\top)^\top,\\
 (P^j,Q^j,\mathcal C)&=\theta(y^j)=(P(y^j)^\top,Q(y^j)^\top, \mathcal C(y^j)^\top)^\top.
 \end{align}
Thus we have $\mathcal C(\varphi_h(y^j))=\mathcal C(y^j)\equiv \mathcal C$ $(j=0,1,\cdots.)$, since the `$\mathcal C$' part in the scheme $\psi_h$ is invariant. Therefore, the scheme $\varphi_h(y)$ preserves the Casimir functions. Next we show it also preserves the Poisson structure, i.e.,
\begin{equation}\label{genstr}
\frac{\partial \varphi_h(y^j)}{\partial y^j}B(y^j)\frac{\partial \varphi_h(y^j)}{\partial y^j}^\top=B(\varphi_h(y^j)), \quad j=0,1,\cdots.
\end{equation}
Denote $\theta(y^j)=:\bar y^j$, and $LHS:=\frac{\partial \varphi_h(y^j)}{\partial y^j}B(y^j)\frac{\partial \varphi_h(y^j)}{\partial y^j}^\top$. Then we have
\begin{align*}
LHS&=\frac{\partial \varphi_h(y^j)}{\partial \psi_h(\bar y^j)}\frac{\partial \psi_h(\bar y^j)}{\partial \bar y^j}\frac{\partial \bar y^j}{\partial  y^j}B(y^j)\frac{\partial \bar y^j}{\partial  y^j}^\top\frac{\partial \psi_h(\bar y^j)}{\partial \bar y^j}^\top \frac{\partial \varphi_h(y^j)}{\partial \psi_h(\bar y^j)}^\top.
\end{align*}
Due to $A(y)B(y)A(y)^\top=B_0$, it holds
\begin{align*}
LHS=\frac{\partial \varphi_h(y^j)}{\partial \psi_h(\bar y^j)}\frac{\partial \psi_h(\bar y^j)}{\partial \bar y^j}\left[\begin{matrix}J^{-1}&0\\0&0\end{matrix}\right]\frac{\partial \psi_h(\bar y^j)}{\partial \bar y^j}^\top \frac{\partial \varphi_h(y^j)}{\partial \psi_h(\bar y^j)}^\top.
\end{align*}
Since $\psi_h(\bar y^j)$ is a symplectic scheme, we have $$\frac{\partial \psi_h(\bar y^j)}{\partial \bar y^j}\left[\begin{matrix}J^{-1}&0\\0&0\end{matrix}\right]\frac{\partial \psi_h(\bar y^j)}{\partial \bar y^j}^\top=\left[\begin{matrix}J^{-1}&0\\0&0\end{matrix}\right],$$ wherefore,
\begin{align*}
LHS=\frac{\partial \varphi_h(y^j)}{\partial \psi_h(\bar y^j)}\left[\begin{matrix}J^{-1}&0\\0&0\end{matrix}\right] \frac{\partial \varphi_h(y^j)}{\partial \psi_h(\bar y^j)}^\top.
\end{align*}
We know that $\psi_h(\bar y^j)=\theta (\varphi_h(y^j))$, thus $\frac{\partial \varphi_h(y^j)}{\partial \psi_h(\bar y^j)}=[A(\varphi_h(y^j))]^{-1}$. Again due to $$A(\varphi_h(y^j))B(\varphi_h(y^j))A(\varphi_h(y^j))^\top=B_0,$$ we have 
\begin{align*}
LHS=B(\varphi_h(y^j)),\quad j=0,1,\cdots.   
\end{align*}
Note that all the derivations above are under `almost surely' sense. \hfill $\square$


\subsection{The $\alpha$-generating function approach for symplectic integration} \label{gen}
The generalized Hamiltonian system \eqref{GSHS} can also be written as \eqref{consps1a}-\eqref{consps1b}, where we only need to solve \eqref{consps1a}. Given initial values $(p,q)$, \eqref{consps1a} can be written as the following standard SHS
\begin{equation}
\begin{split}
\label{SHS}
&dP=-\frac{\partial H_0(P,Q)}{\partial Q}dt-\sum\limits_{r=1}^{m}\frac{\partial H_r(P,Q)}{\partial Q}\circ dW_r(t),\quad P(t_0)=p,\\
&dQ=\frac{\partial H_0(P,Q)}{\partial P}dt+\sum\limits_{r=1}^{m}\frac{\partial H_r(P,Q)}{\partial P}\circ dW_r(t),\quad Q(t_0)=q,
\end{split}
\end{equation}
where $P,Q,p,q\in \R^n$. 
Assume that the Hamiltonian functions $H_{r}$ ($r=0, \dots, m$) belong to $C^{\infty}.$ 
In addition, we also suppose that, for any $(P,Q)\in\R^{2n}$, $(\bar{P},\bar{Q})\in \R^{2n}$, there exist $L_1>0$ and $L_2>0$ such that
\[
\begin{aligned}
&\sum_{r=0}^{m}\left(\left|\nabla_{p} H_{r}(P, Q)-\nabla_{p} H_{r}(\bar{P}, \bar{Q})\right|+\left|\nabla_{q} H_{r}(P, Q)-\nabla_{q} H_{r}(\bar{P}, \bar{Q})\right|\right)\\
+&\frac{1}{2}\sum_{r=1}^{m} \left| \sigma_{r}'(P,Q) \sigma_{r}(P,Q)- \sigma_{r}'(\bar{P},\bar{Q}) \sigma_{r}(\bar{P},\bar{Q})\right|
\leq  L_{1}(|P-\bar{P}|+|Q-\bar{Q}|),
\end{aligned}
\]
and
\[
\sum_{r=0}^{m}\left(\left|\nabla_{p} H_{r}(P, Q)\right|+\left|\nabla_{q} H_{r}(P, Q)\right|\right)+\frac{1}{2}\sum_{r=1}^m |\sigma_r'(P,Q)\sigma_r(P,Q)| \leq L_{2}(1+|P|+|Q|),
\]
where $\sigma_{r}(P,Q)=\left(-{\nabla_{q} H_{r}(P, Q)^\top,\nabla_{p} H_{r}(P, Q)^\top}\right)^\top$.
The above two conditions guarantee the local existence and uniqueness of the solution of the SHS \eqref{SHS}.

The phase flow of \eqref{SHS} preserves the symplectic structure (\cite{milbook,mil1,mil2}), which, 
using the differential 2-form, can be characterized as
\begin{equation*}
dP(t)\wedge dQ(t)=dp\wedge dq,\quad \forall t\geq t_0.
\end{equation*}
A symplectic numerical method $\{P_k,Q_k\}_{k}$ with $(P_0,Q_0)=(p,q)$ is a method that can preserve the symplectic structure, namely, 
\begin{equation}\label{symmethod}
dP_{k+1}\wedge d Q_{k+1}=dP_k\wedge d Q_k, \quad \forall \,\,\,k\ge 0.
\end{equation}
As was shown in \cite{fengq,ECG2002,wang,deng}, a mapping $(p^\top,q^\top)^\top\rightarrow (P^\top,Q^\top)^\top$ is symplectic if there exists a locally smooth generating function $S(q,Q,t)$, such that 
\begin{align}
\label{todiff}
P^\top dQ-p^\top dq=dS(q,Q)
\end{align}
for every fixed $t$. In stochastic case, the generating function $S(q,Q,t,\omega)$ can be obtained by solving the stochastic Hamilton-Jacobi partial differential equation \cite{bismut,deng,wang,wangh} 
\begin{align}\label{shj}
\partial_{t}S(q,Q,t,\omega)=
-H_0\left(\frac{\partial S}{\partial Q},Q\right)dt
-\sum_{r=1}^{m}H_r\left(\frac{\partial S}{\partial Q},Q\right)\circ dW_r(t),
\end{align}
with initial conditions $\frac{\partial S}{\partial Q_i}(q,q,t_0)+\frac{\partial S}{\partial q_i}(q,q,t_0)=0$ $(i=1,\cdots,n)$. The notion $S(q, Q,t,\omega)$ represents a family of real valued stochastic processes with parameter $(q,Q)\in \R^{2n}$, which can be regarded as a random field with parameters $(q,Q,t)$ (\cite{deng}). If $S(q,Q,t,\omega)$ is a $C^{\infty} $function of $(q,Q)$ for almost every $\omega$ for each $t$, it can be regarded as a $C^{\infty}$ value process (\cite{kunita1997}). It can be proved that, under certain conditions( see \cite{bismut,kunita1997,deng}), a local solution $S(q,Q,t,\omega)$ of \eqref{shj} can almost surely generate the flow $\varphi_t: (p^\top,q^\top)^\top\rightarrow (P(t)^\top,Q(t)^\top)^\top$ $(t\in[t_0,\tau])$ of the SHS \eqref{SHS} via the relation   
\begin{align*}
P(t)=\frac{\partial S(q,Q(t),t)}{\partial Q},
\quad p=-\frac{\partial S(q,Q(t),t)}{\partial q},
\end{align*}
if the matrix $\left({\partial ^2S}\over{\partial q_i\partial Q_j}\right)$ is almost surely invertible in $t\in[t_0,\tau]$ where $\tau>t_0$ is a stopping time.

In addition to the aforementioned generating function $S(q,Q,t)$, with different coordinates, there can be other kinds of generating functions (\cite{fengq,ECG2002,wang,deng}). We unify and extend them to the $\alpha$-generating functions with parameter $\alpha\in\left[0,1\right]$ in the following.

Denote $\hat{P}=(1-\alpha)p+\alpha P,$ $\hat{Q}=(1-\alpha)Q+\alpha q$ with $\alpha\in\left[0,1\right].$  We have the following theorem regarding the $\alpha$-generating function $\hat S_{\alpha}(\hat{P}, \hat{Q},t)$.

\begin{tm}
\label{HATGF}
A mapping $(p^\top,q^\top)^\top\rightarrow (P^\top,Q^\top)^\top$ is symplectic if there exist the generating functions $\hat{S}_{\alpha}(\hat{P}, \hat{Q},t)$ $(\alpha\in \left[0,1\right])$, such that the following equations hold for every fixed $t$, 
\begin{align}
 a)&p^\top d\hat{Q}+Q^\top d\hat{P}=d\left[\hat{P}^\top\hat{Q}\right]+\alpha d\hat{S}_{\alpha},\quad\mbox{if}\,\,\alpha\in(0,1];\label{1}\\
b)& P^\top d\hat Q+q^\top d\hat P=d\left[\hat P^\top\hat Q\right]-(1-\alpha)d \hat S_{\alpha},\quad \mbox{if}\,\,\alpha\in [0,1); \label{2}\mbox{or uniformly}\\
c)& \left[{\bf 1}_{(0,1]}(\alpha)p^\top+{\bf 1}_{[0,1)}(\alpha)P^\top\right]d\hat Q+\left[{\bf 1}_{(0,1]}(\alpha)Q^\top+{\bf 1}_{[0,1)}(\alpha)q^\top\right]d\hat P \nonumber\\&=\left[{\bf 1}_{(0,1]}(\alpha)+{\bf 1}_{[0,1)}(\alpha)\right]d\left[\hat P^\top \hat Q\right]+\left[\alpha {\bf 1}_{(0,1]}(\alpha)-(1-\alpha){\bf 1}_{[0,1)}(\alpha)\right] d\hat S_{\alpha}, \nonumber \\& \hspace*{8.1cm} \quad \mbox{for}\,\,\,\alpha\in[0,1]. \label{3}
\end{align}
\end{tm}
{\bf Proof.}
 ${\bf 1}_{A}(\alpha)$ denotes the indicator function of the set $A$. We first consider the case for $\alpha\in(0,1).$ Note that \eqref{todiff} is valid (see e.g. \cite{ECG2002,wang,deng}). 

Multiplying both sides of \eqref{todiff} by $\alpha(1-\alpha),$ we obtain
$$\alpha P^\top d\left[(1-\alpha)Q\right]-(1-\alpha)p^\top d\left[\alpha q\right]=\alpha(1-\alpha)dS.$$
Adding the term $\alpha P^\top d\left[\alpha q\right]-(1-\alpha)p^\top d\left[(1-\alpha)Q\right]$ to both sides of the equation above, we get
\begin{align*}
\alpha P^\top d\hat{Q}-(1-\alpha)p^\top d\hat{Q}
=\alpha(1-\alpha)dS
+\alpha P^\top d\left[\alpha q\right]
-(1-\alpha)p^\top d\left[(1-\alpha)Q\right],
\end{align*}
which leads to
\begin{align*}
&d\hat{Q}^\top(\alpha P)=\alpha P^\top d\hat{Q}+\hat{Q}^\top d\left[\alpha P\right]\\
=&\alpha(1-\alpha)dS+\alpha P^\top d\left[\alpha q\right]-(1-\alpha)p^\top d\left[(1-\alpha)Q\right]+(1-\alpha)p^\top d\hat{Q}\\	
&
+(1-\alpha)Q^\top d\hat{P}-(1-\alpha)Q^\top d\left[(1-\alpha)p\right]+\alpha q^\top d\left[\alpha P\right]\\
=&\alpha(1-\alpha)dS+(1-\alpha)(p^\top d\hat Q+Q^\top d\hat P)
-(1-\alpha)^2d\left[Q^\top p\right]+\alpha^2d\left[P^\top q\right].
\end{align*}
Based on the fact that
$$
\alpha P^\top d\hat{Q}+\hat{Q}^\top d\left[\alpha P\right]
+(1-\alpha)^2d\left[Q^\top p\right]-\alpha^2d\left[P^\top q\right]
=(1-\alpha)d\left[\hat P^\top Q\right],
$$
we have
\begin{equation}\label{two}
\alpha(1-\alpha)dS+(1-\alpha)(p^\top d\hat Q+Q^\top d\hat P)=(1-\alpha)d\left[\hat P^\top Q\right].
\end{equation}
For $\alpha\in (0,1)$, it can be derived that \eqref{two} is equivalent to both of the following equations
\begin{align}
&p^\top d\hat Q+Q^\top d\hat P=d \left[\hat P^\top Q\right]-\alpha dS,\label{two1}\\
&P^\top d \hat Q+q^\top d\hat P=d\left[\hat P^\top q\right]+(1-\alpha)dS,\label{two2}
\end{align}
via eliminating $1-\alpha$ or $\alpha$ from \eqref{two}, respectively.
\eqref{two1} implies that there exists function 
\begin{equation}\label{s1}
\hat S_{\alpha}=\hat P^\top(Q-q)-S,
\end{equation}
such that \eqref{1} holds for $\alpha\in (0,1)$, namely,
\begin{align*}
p^\top d\hat Q+Q^\top d\hat P=d\left[\hat P^\top \hat Q\right]+\alpha d\hat S_{\alpha},
\end{align*}
and \eqref{two2} suggests to let 
\begin{equation}\label{s2}
(1-\alpha)\hat S_{\alpha}=\hat P^\top (\hat Q-q)-(1-\alpha)S,
\end{equation}
which then satisfies \eqref{2} for $\alpha\in (0,1)$, i.e.,
$$P^\top d\hat Q+q^\top d\hat P=d(\hat P^\top \hat Q)-(1-\alpha)d \hat S_{\alpha}.$$

It is not difficult to see from the derivation that, for $\alpha\in(0,1)$, \eqref{1} and \eqref{2} are equivalent, and \eqref{s1} and \eqref{s2} are equivalent as well. If $\alpha=1$, we can check that $\hat S_1(P,q,t)$ is just the first kind of generating function $S^1(P,q,t)$ (\cite{ECG2002,wang,deng}), which satisfies \eqref{1} and \eqref{s1}. If $\alpha=0$, $\hat S_0(p,Q,t)$ corresponds to the second kind of generating function $S^2(p,Q,t)$ (\cite{ECG2002,Anton2}), which satisfies \eqref{2} and \eqref{s2}. Thus \eqref{1} and \eqref{2} are proved, and \eqref{3} is a naturally unified expression of \eqref{1} and \eqref{2} for all $\alpha\in[0,1]$. 
\hfill$\square$
\begin{rk}
When $\alpha={{1}\over{2}}$, the function $\hat S_{{1}\over{2}}(\hat P,\hat Q,t)$ is the third kind of generating function \\$S^3(\frac{P+p}{2}, \frac{Q+q}{2},t)$ (see e.g. \cite{ECG2002}).
\end{rk}

Similar to $S(q,Q,t,\omega),$ the generating functions $\hat S_{\alpha}(\hat P,\hat Q,t,\omega)$ with $\alpha\in\left[0,1\right]$ can also be associated with the stochastic Hamilton-Jacobi partial differential equation.
Following a similar procedure of proving Theorem 2.1 in \cite{deng}, we can prove the following theorem.

\begin{tm}
Let $\hat S_{\alpha}(\hat P,\hat Q,t,\omega)$ ($\alpha\in[0,1]$) be a locally smooth solution of the stochastic Hamilton-Jacobi partial differential equation
\begin{align}\label{shjpde}
\partial_{t}\hat S_{\alpha}(\hat P, \hat Q,t,\omega)
=
\sum_{r=0}^{m}H_r\left(\hat P-(1-\alpha)\frac{\partial \hat S_{\alpha}}{\partial \hat Q},
\hat{Q}+\alpha\frac{\partial \hat S_{\alpha}}{\partial \hat P}\right)\circ dW_r(t)
\end{align}
with initial value $\hat S_{\alpha}(\hat P,\hat Q,t_0)=0$, $dW_0(t)=dt,$ such that almost sure $\hat S_{\alpha}(\hat P,\hat Q,t,\omega), \partial \hat S_{\alpha}(\hat P,\hat Q,t,\omega) / \partial \hat P$ and $\partial \hat S_{\alpha}(\hat P,\hat Q,t,\omega)/ \partial \hat Q$ are local Stratonovich semi-martingales, continuous on $(\hat{P}, \hat{Q}, t)$ and $C^{\infty}$ value processes. If in addition there exists a stopping time $\tau>t_0$ almost surely such that the matrix $\left(\frac{\partial^2 \hat S_{\alpha}(\hat P,\hat Q,t,\omega)}{\partial \hat P_i\partial \hat{Q}_j}\right)$
is almost surely invertible for $t_0\leq t<\tau,$
then the mapping $(p,q)\mapsto(P(t,\omega),Q(t,\omega))$ $(t_0\leq t< \tau)$ defined by
\begin{align}\label{relation}
P(t,\omega)=p-\frac{\partial \hat S_{\alpha}(\hat P,\hat Q,t,\omega)}{\partial \hat{Q}},
\quad Q(t,\omega)=q+\frac{\partial \hat S_{\alpha}(\hat P,\hat Q,t,\omega)}{\partial \hat{P}}
\end{align}
is the flow of the SHS \eqref{SHS}.
\end{tm}

The integral form of the stochastic Hamilton-Jacobi PDE (\ref{shjpde}) under its initial condition is
\begin{align}\label{intshj}
\hat{S}_{\alpha}(\hat{P},\hat{Q},t,\omega)=\sum_{r=0}^m\int_{t_0}^t H_r\left(\hat P-(1-\alpha)\frac{\partial \hat S_{\alpha}}{\partial \hat Q},
\hat{Q}+\alpha\frac{\partial \hat S_{\alpha}}{\partial \hat P}\right)\circ dW_r(s),
\end{align}
where $\hat{P},\hat{Q}$ are regarded as parameters. Following the idea for deterministic case in \cite{fengq}, since $H_r$ ($r=0,\dots,m$) are assumed to be $C^{\infty}$, we can perform a Stratonovich-Taylor expansion of (\ref{intshj}) by expanding the integrands $H_r\left(\hat P-(1-\alpha)\frac{\partial \hat S_{\alpha}}{\partial \hat Q},\hat{Q}+\alpha\frac{\partial \hat S_{\alpha}}{\partial \hat P}\right)$ at $(\hat{P},\hat{Q})$, which will assume the following formal series expansion of $\hat{S}_{\alpha}$:

\begin{align}
\label{GF}
\hat S_{\alpha}(\hat P,\hat Q,t,\omega)=\sum_{\gamma} G_{\gamma}^{\alpha}(\hat P,\hat Q)J_{\gamma},
\end{align}
where 
\begin{align}
\label{int;ssm}
J_{\gamma}=\int_0^t\int_0^{s_l}\cdots\int_0^{s_2}\circ dW_{j_1}(s_1)\circ dW_{j_2}(s_2)\circ\cdots\circ dW_{j_l}(s_l)
\end{align}
with multi-index $\gamma=(j_1,j_2,\cdots,j_l)$, $j_i\in\{0,1,\cdots,m\}$, $(i=1,\cdots,l),$ $l\geq 1.$
To determine the coefficients $G_{\gamma}^{\alpha}(\hat P,\hat Q)$ in \eqref{GF}, one can substitute the ansatz (\ref{GF}) into (\ref{intshj}) to compare like powers of $t$. To this end, we first introduce the following notations:
\begin{itemize}
\item Denote by $l(\gamma)$ and $\gamma-$ the length of $\gamma$ and the multi-index resulted from discarding the last index of $\gamma,$ respectively.
\item Define $\gamma\ast \gamma'=(j_1,\cdots,j_l,j_1',\cdots,j_{l'}')$ where $\gamma=(j_1,\cdots,j_l)$ and $\gamma'=(j_1',\cdots,j_{l'}').$
\item Let
\begin{equation*}
\Lambda_{ \gamma_1, \gamma_2}=\left\{
\begin{split}
\{(j_1,j_1'),(j_1',j_1)\},\quad &\mbox{if}\quad l=l'=1,\\
\{\Lambda _{(j_1),\gamma_2- }*(j_{l'}'),\gamma_2\ast (j_1)\},\quad &\mbox{if}\quad l=1, l'\neq 1,\\
\{\Lambda _{\gamma_1-,(j_1')}\ast(j_{l}),\gamma_1\ast (j_{1}')\},\quad &\mbox{if}\quad l\neq1, l'=1,\\
\{\Lambda_{\gamma_1-,\gamma_2}\ast (j_l),\Lambda_{\gamma_1,\gamma_2-}\ast (j_{l'}')\},\quad &\mbox{if}\quad l\neq 1, l'\neq 1,\end{split}\right.
\end{equation*}
where the concatenation $'\ast'$ between a set of multi-indices $\Lambda$ and $\gamma$ is $\Lambda\ast\gamma=\{\beta*\gamma|\beta\in \Lambda\}.$
\item For $k>2$, $\Lambda_{\gamma_1,\cdots, \gamma_k}=\{\Lambda_{\beta, \gamma_k}|\beta\in\Lambda_{\gamma_1,\cdots, \gamma_{k-1}}\}$.
\end{itemize}

Now, we use the same technique in \cite{deng}, to substitute the series expansion \eqref{GF} into the stochastic Hamilton-Jacobi partial differential equation \eqref{shjpde} and take Taylor's series expansions of $H_r$ at $(\hat P,\hat Q)$ $(r=0,\cdots,m)$, to obtain the following expression of $\hat S_{\alpha}(\hat{P},\hat{Q},t)$:
\begin{align}\label{sal}
&\hat S_{\alpha}(\hat{P},\hat{Q},t)= \nonumber \\
&\sum_{r=0}^{m}\int_0^t H_r(\hat{P},\hat{Q})\circ dW_r(s)+\sum\limits_{i=1}^{\infty}\frac{1}{i!}\sum_{k_1,\cdots,k_i=1}^{d}\sum_{j=0}^{i}\frac{\partial^i H_r(\hat{P},\hat{Q})}{\partial \hat{P}_{k_1}\cdots\partial \hat{P}_{k_j}\partial \hat{Q}_{k_{j+1}}\cdots\partial \hat{Q}_{k_i}} \nonumber\\
&\quad\cdot 1C_i^j(\alpha-1)^j\alpha^{i-j}\sum_{\gamma_1,\cdots,\gamma_i}
\frac{\partial G_{\gamma_1}}{\partial \hat Q_{k_1}}
\cdots\frac{\partial G_{\gamma_j}}{\partial \hat Q_{k_j}}
\frac{\partial G_{\gamma_{j+1}}}{\partial \hat P_{k_{j+1}}}
\cdots\frac{\partial G_{\gamma_i}}{\partial \hat P_{k_i}}
\int_{0}^{t}\prod_{k=1}^{i}J_{\gamma_k}\circ dW_r(s).
\end{align}
Due to the relation (\cite{kloeden})
\begin{align*}
\prod_{k=1}^i J_{\gamma_k}=\sum_{\beta\in\Lambda _{\gamma_1,\cdots,\gamma_n}}J_{\beta},
\end{align*}
and after equating coefficients on both sides of the equation \eqref{sal}, we obtain
\begin{equation}
\begin{split}
\label{2.16}
G_{\gamma}^{\alpha}=&\sum_{i=1}^{l(\gamma)-1}\frac{1}{i!}\sum_{k_1,\cdots,k_i=1}^{d}\sum_{j=0}^{i}\frac{\partial^i H_r(\hat{P},\hat{Q})}{\partial \hat{P}_{k_1}\cdots\partial \hat{P}_{k_j}\partial \hat{Q}_{k_{j+1}}\cdots\partial \hat{Q}_{k_i}}C_i^j(\alpha-1)^j\alpha^{i-j}\\
&\cdot
{\scriptsize
\sum_{\begin{array}{c}l(\gamma_1)+\cdots+l(\gamma_i)=l(\gamma)-1\\ \gamma-\in\Lambda \gamma_1,\cdots,\gamma_i\end{array}}
}
\frac{\partial G_{\gamma_1}}{\partial \hat Q_{k_1}}
\cdots\frac{\partial G_{\gamma_j}}{\partial \hat Q_{k_j}}
\frac{\partial G_{\gamma_{j+1}}}{\partial \hat P_{k_{j+1}}}
\cdots\frac{\partial G_{\gamma_i}}{\partial \hat P_{k_i}}
\end{split}
\end{equation}
for $\gamma=(i_1,\cdots,i_{l-1},r)$ with $l>1$, $i_1,\cdots,i_{l-1},r$ taking values from $\{0,1,\cdots,m\}$ without duplication.
If there are duplicates in $\gamma$, one can still use the formula after assigning different subscripts to the duplicates.
For $l(\gamma)=1$, i.e., $\gamma=(r)$, $ G_r^{\alpha}=H_r(\hat P, \hat Q).$
In sum, the generating function $\hat S_{\alpha}$ can be expressed as
\begin{align*}
\hat S_{\alpha}
=&H_0J_0+\sum\limits_{r=1}^{m}H_rJ_r
+(2\alpha-1)\sum\limits_{k=1}^{n}
\frac{\partial H_0}{\partial \hat Q_k}
\frac{\partial H_0}{\partial \hat P_k}J_{(0,0)}\\
&+\sum\limits_{r=1}^{m}\sum\limits_{s=1}^{m}\sum\limits_{k=1}^{n}
\left(
\alpha\frac{\partial H_r}{\partial \hat Q_k}
\frac{\partial H_s}{\partial \hat P_k}
+(\alpha-1)\frac{\partial H_r}{\partial \hat P_k}
\frac{\partial H_s}{\partial \hat Q_k}
\right)J_{(s,r)}\\
&+\sum\limits_{r=1}^{m}\sum\limits_{k=1}^{n}
\left(
\alpha\frac{\partial H_r}{\partial \hat Q_k}
\frac{\partial H_0}{\partial \hat P_k}
+(\alpha-1)\frac{\partial H_r}{\partial \hat P_k}
\frac{\partial H_0}{\partial \hat Q_k}
\right)J_{(0,r)}\\
&+\sum\limits_{r=1}^{m}\sum\limits_{k=1}^{n}
\left(
(\alpha-1)\frac{\partial H_0}{\partial \hat P_k}
\frac{\partial H_r}{\partial \hat Q_k}
+\alpha\frac{\partial H_0}{\partial \hat Q_k}
\frac{\partial H_r}{\partial \hat Q_k}
\right)J_{(r,0)}+\cdots.
\end{align*}

To construct a symplectic numerical scheme with desired mean-square order $\mathcal K$ via truncating the generating functions, \cite{Anton2,deng} proposed the following procedure: First replace every multiple Stratonovich integral in a generating function $S$ by its equivalent combination of multiple It\^o integrals.
Then, truncate the series of $S$ to include all terms containing It\^o integrals with multi-index $\gamma$ belonging to the set $\AAA_{\mathcal K}:=\{\gamma:l(\gamma)+n(\gamma)\leq 2\mathcal K,{\rm{or}},l(\gamma)=n(\gamma)=\mathcal K+0.5\}$
with $n(\gamma)$ being the number of zero components in $\gamma.$

Regarding our $\alpha$-generating function $\hat S_{\alpha}$, for example, if $\mathcal K=1$ and $m=1,$ the truncated generating function is
\begin{align*}
\bar S_{\alpha}=&H_0I_{(0)}+H_1I_{(1)}+(2\alpha-1)\frac{\partial H_1}{\partial \hat Q}
\frac{\partial H_1}{\partial \hat P}\left(I_{(1,1)}+\frac 12 I_{(0)}\right),
\end{align*}
which produces the following symplectic schemes according to the relation \eqref{relation}
\begin{align}
\label{sns1}
\begin{bmatrix}
P_{n+1}\\
Q_{n+1}
\end{bmatrix}=
\begin{bmatrix}
P_n\\
Q_n
\end{bmatrix}
+J^{-1}\nabla
\bar S((1-\alpha)P_n+\alpha P_{n+1},(1-\alpha)Q_{n+1}+\alpha Q_{n}).
\end{align}

In fact, \eqref{sns1} is the same as the $\theta$-method introduced in \cite{milbook}, where its mean-square convergence rate is given based on the fundamental theorem on mean-square convergence. 

The $\alpha$-generating function approach enriches the generating function theory of constructing symplectic schemes for Hamiltonian systems. It allows continuously varying choice of $\alpha$ from $[0,1]$, and creates a large class of symplectic integrators.  The $\alpha$-generating function approach itself is of theoretical and practical significance, though we only embed it in this paper into the integration strategy for SPSs, to construct symplectic methods for the SHSs resulted from the canonical transformation acted on the SPSs. 

\subsection{Applications to the SRB and SLV systems}
\label{app}
We use our integration strategy for stochastic Poisson systems to solve numerically a stochastic rigid body system and a stochastic Lotka-Volterra system. We first set up appropriate coordinate transformation to transform the SPSs to their canonical forms, i.e., the generalized stochastic Hamiltonian systems, and use the $\alpha$-generating function method to create symplectic schemes for the SHSs. Then we perform the inverse coordinate transformation on the symplectic schemes to get the Poisson integrators for the original SPSs, which we call the `$\alpha$-generating schemes' for brevity. Certain non-canonical coordinate transformation method will also be illustrated.

\subsubsection{The three-dimensional stochastic rigid body system}\label{apprb}
Recall the stochastic rigid body system 
\begin{align}
\label{model;rigid}
dy=B_1(y)\nabla K_1(y)(dt+c_1\circ dW(t)),\quad y(0)=(y_1^0,y_2^0,y_3^0)^\top,
\end{align}
where
$K_1=\frac 12
\left(
\frac{y_1^2}{I_1}
+\frac{y_2^2}{I_2}
+\frac{y_3^2}{I_3}
\right)$, $I_1,I_2,I_3,c_1$ are constants, and
$$B_1=
\begin{bmatrix}
0& -y_3 & y_2\\
y_3 & 0& -y_1\\
-y_2 & y_1& 0
\end{bmatrix}.$$
It possesses the Casimir function $$C_1(y)={{1}\over{2}}(y_1^2+y_2^2+y_3^2)\equiv {{1}\over {2}}\left((y_1^0)^2+(y_2^0)^2+(y_3^0)^2\right)=:\mathcal C_1.$$

First we look for a coordinate transformation $\bar{y}(y)=(\bar y_1(y),\bar y_2(y),\bar y_3(y))^\top$ with  invertible Jacobian matrix $$A(y) = \left(A_{i j}(y)\right)=\left({\partial \bar y_i}\over{\partial y_j}\right)\quad (i,j=1,2,3)$$ such that 
\begin{equation}\label{catr}
A(y)B_1(y) A(y)^{\top}=B_0,
\end{equation}
where $B_0$ can be
\begin{align*}
\begin{bmatrix} 0 & -1 & 0 \\ 1 & 0 & 0 \\ 0 & 0 & 0 \end{bmatrix},\quad \quad
\begin{bmatrix} 0 & 0 & -1 \\ 0 & 0 & 0 \\ 1 & 0 & 0 \end{bmatrix},\quad {\rm or}\quad
\begin{bmatrix} 0 & 0 & 0 \\ 0 & 0 & -1 \\ 0 & 1 & 0 \end{bmatrix}.
\end{align*}
Now we take the first matrix above to be $B_0$.  Then \eqref{catr} is equivalent to the following equations with respect to the Poisson bracket defined by $B_1$
\begin{equation}\label{e3}
\begin{split}
&\{\bar y_1,\bar y_1\}=0,\quad \{\bar y_1,\bar y_2\}=-1, \quad \{\bar y_1,\bar y_3\}=0,\\
&\{\bar y_2,\bar y_1\}=1,\quad \{\bar y_2,\bar y_2\}=0, \quad \{\bar y_2,\bar y_3\}=0,\\
&\{\bar y_3,\bar y_1\}=0,\quad \{\bar y_3,\bar y_2\}=0, \quad \{\bar y_3,\bar y_3\}=0.
\end{split}
\end{equation}
Due to anti-symmetry of the Poisson bracket, the nine equations above can be reduced to the following three equations
\begin{equation}\label{eq3}
\{\bar y_2,\bar y_1\}=1,\quad \{\bar y_3,\bar y_1\}=0,\quad \{\bar y_3,\bar y_2\}=0.
\end{equation}
The last two equations above imply that we can choose $\bar y_3=\mathcal C_1$, according to the property of the Casimir functions. The first equation can be expressed explicitly as
\begin{equation}\label{eqq1}
(A_{12} A_{23} - A_{13} A_{22}) y_1 + (A_{13} A_{21} - A_{11} A_{23}) y_2 + (A_{11} A_{22} - A_{12} A_{21}) y_3 =1.
\end{equation}
This is actually a partial differential equation with respect to $\bar y_1(y_1,y_2,y_3)$ and $\bar y_2(y_1,y_2,y_3)$, which possesses possibly many variants of solutions. 
If we let, e.g.,  $\bar{y}_1 = y_2$, then the equation \eqref{eqq1} becomes
\begin{equation}\label{eqq2}
A_{23} y_1 - A_{21} y_3 = 1,
\end{equation}
and it can be verified that
$\bar{y}_2 = \arctan \left(\frac{y_3}{y_1}\right)$ solves the equation \eqref{eqq2}.
Thus, we find the following coordinate transformation
\begin{equation}\label{trans1}
\bar{y}_1 = y_2, \quad
\bar{y}_2 = \arctan \left(\frac{y_3}{y_1}\right), \quad \bar{y}_3 = \mathcal C_1,
\end{equation}
and its inverse 
\begin{equation}\label{intrans2}
y_1 = \sqrt{2 \mathcal C_1 - \bar{y}_1^2} \cos (\bar{y}_2), \quad y_2 = \bar{y}_1, \quad y_3 = \sqrt{2 \mathcal C_1 - \bar{y}_1^2} \sin (\bar{y}_2).
\end{equation}
Simultaneously, we obtain the stochastic Hamiltonian system of $\bar y_1$ and $\bar y_2$
\begin{align}
\label{model:ESHS1}
d\begin{bmatrix}
\bar y_1\\
\bar y_2
\end{bmatrix}=
\begin{bmatrix}
0& -1\\
1& 0
\end{bmatrix}
\nabla H(\bar y_1,\bar y_2)(dt+c_1\circ dW(t)),
\end{align}
where $$H(\bar y_1,\bar y_2)=
\frac{1}{2I_1} (2\mathcal C_1-\bar y_1^2)\cos^2 (\bar y_2)
+\frac{1}{2I_2} \bar y_1^2
+\frac{1}{2I_3} (2\mathcal C_1-\bar y_1^2)\sin^2 (\bar y_2).$$

Next we apply the symplectic scheme \eqref{sns1}, which is given by the $\alpha$-generating function approach and of mean-square order 1, to the SHS \eqref{model:ESHS1}. Substituting the derivatives of the Hamiltonian function
\begin{align*}
&\frac{\partial H}{\partial \bar y_1}
=\left(\frac{1}{I_2}-
\frac{\cos^2 (\bar y_2)}{I_1}
-\frac{\sin^2 (\bar y_2)}{I_3}\right)\bar y_1,\quad \frac{\partial H}{\partial \bar y_2}
=\left(\frac{1}{2I_3}-
\frac{1}{2I_1}\right)(2L-\bar y_1^2)\sin (2\bar y_2),\\
&\frac{\partial^2 H}{\partial \bar y_1^2}
=\frac{1}{I_2}-
\frac{\cos^2 (\bar y_2)}{I_1}
-\frac{\sin^2 (\bar y_2)}{I_3},\qquad\,\,\,\,\,\,
\frac{\partial^2 H}{\partial \bar y_2^2}
=\left(\frac{1}{I_3}-
\frac{1}{I_1}\right)
(2L-\bar y_1^2)\cos (2\bar y_2),\\
&\frac{\partial^2 H}{\partial \bar y_1\partial \bar y_2}
=\left(\frac{1}{I_1}-
\frac{1}{I_3}\right)\bar y_1\sin (2\bar y_2)
\end{align*}
into \eqref{sns1}, we obtain the following symplectic scheme
\begin{equation}
\label{RGnm1}
\begin{split}
P_{n+1}=&P_n-
\left(
\frac{1}{2I_3}-
\frac{1}{2I_1}
\right)
\left(2\mathcal C_1-\bar P_n^2\right)
\sin \left(2\bar Q_n\right)
(h+c_1\Delta W_n)\\
&+c_1^2C_{\alpha}
(2\mathcal C_1-\bar P_n^2)\bar P_n\left(
\cos (2\bar Q_n)
\left(
\frac{1}{I_2}-
\frac{\cos^2 (\bar Q_n)}{I_1}
-\frac{\sin^2 (\bar Q_n)}{I_3}
\right)\right.\\
&\left.-\sin ^2(2\bar Q_n)\left(\frac{1}{2I_3}-\frac{1}{2I_1}\right)\right)
\Delta W_n^2
\\
Q_{n+1}
=&Q_n+\left(\frac{1}{I_2}-
\frac{\cos^2 (\bar Q_n)}{I_1}
-\frac{\sin^2 (\bar Q_n)}{I_3}\right)\bar P_n(h+c_1\Delta W_n)\\
&+c_1^2C_\alpha
\left(
\frac{1}{I_2}-
\frac{\cos^2 (\bar Q_n)}{I_1}
-\frac{\sin^2 (\bar Q_n)}{I_3}
\right)
\left(\frac{3}{2}\bar P_n^2-\mathcal C_1\right)\sin (2\bar Q_n)
\Delta W_n^2
\end{split}
\end{equation}
where $\bar P_n=(1-\alpha)P_n+\alpha P_{n+1},$ $\bar Q_n=\alpha Q_n+(1-\alpha) Q_{n+1},$ $C_\alpha=\left(\alpha-\frac12\right)\left(
\frac 1{I_1}-\frac 1{I_3}\right)$ with $\alpha\in\left[0,1\right].$

These are symplectic schemes which are implicit. To fix the problems caused by the unboundedness of $\Delta W_{n}=\sqrt{h}\xi_{n}$, we follow the method given in \cite{mil2} to truncate the $\mathcal N(0,1)$-distributed random variable  $\xi_n$ to another bounded random variable $\zeta_{n}$. In detail,
\begin{equation*}
\label{trun}
\zeta_{n}=\left\{\begin{array}{l}\ \xi_{n},\,\,\,\,\mbox{if}\,\,\,\,|\xi_{n}| \leq A_{h},\\ \ A_{h},\,\,\,\,\mbox{if}\,\,\,\, \xi_{n} > A_{h},\\ \ -A_{h},\,\,\,\,\mbox{if}\,\,\,\, \xi_{n} < -A_{h},\end{array}\right.
\end{equation*}
where $A_{h}=\sqrt{2k|\ln h|}$, $k \geq 1$. It is also indicated in \cite{mil2} that, the truncation error can be merged into the error of the numerical scheme by choosing sufficiently large parameter $k$, which should be at least $2\mathcal K$ if the numerical scheme containing such a truncation is expected to possess root-mean-square convergence order $\mathcal K$. In our numerical tests in Section \ref{nmil} we take $k=4$.

By the inverse coordinate transformation \eqref{intrans2}, we get the following $\alpha$-generating schemes for the original stochastic rigid body system \eqref{model;rigid}
\begin{align}
\label{snm;rg1}
Y^1_n= \sqrt{2 \mathcal C_1 - P_n^2} \cos (Q_n), \quad
Y^2_n= P_n, \quad
Y^3_n = \sqrt{2 \mathcal C_1 - P_n^2} \sin (Q_n).
\end{align}

It is easy to see that \eqref{snm;rg1} preserves the Casimir function, since
$$(Y^1_n)^2+(Y^2_n)^2+(Y^3_n)^2=2\mathcal C_1.$$
Moreover, it inherits the Poisson structure of the stochastic rigid body system \eqref{model;rigid}, according to the proof of Theorem \ref{gen1}.

On the other hand, the quadratic form of the Casimir function 
$$C_1(y)=\frac{1}{2}(y_1^2+y_2^2+y_3^2)\equiv \mathcal C_1$$
motivates a spherical coordinate transformation $\phi:(\theta_1,\theta_2)\rightarrow(y_1,y_2,y_3)$, i.e.
\begin{equation}
\label{sps-eq-9}
y_1 = R \cos \theta_1 \cos \theta_2,
\quad
y_2 = R \cos \theta_1 \sin \theta_2, \quad
y_3 = R\sin \theta_1,
\end{equation}
where $R = \sqrt {2\mathcal C_1}$.

Using the inverse mapping of $\phi$, we have that
\begin{equation}
\begin{split}
\label{sps-eq-10}
d\theta_1=& R \left(\frac{1}{I_2}-\frac{1}{I_1}\right)
\cos\theta_1 \sin\theta_2 \cos\theta_2 (dt+c_1\circ dW(t)), \\
d\theta_2=& R \sin\theta_ 1 \left(
\Big{(}\frac{1}{I_1}-\frac{1}{I_3}\Big{)}
\cos^2 \theta_2
- \Big{(}\frac{1}{I_3}-\frac{1}{I_2}\Big{)}
\sin^2 \theta_2
\right)
(dt+c_1 \circ dW(t)).
\end{split}
\end{equation}
We can apply the midpoint rule, which corresponds to the scheme \eqref{sns1} with $\alpha=\frac{1}{2}$, to the system \eqref{sps-eq-10}, to get
\begin{equation} 
\begin{split}
\label{snm;sym2}
\Theta^1_{n+1}&=\Theta^1_{n}+ R \left(\frac{1}{I_2}-\frac{1}{I_1}\right)
\cos\bar \Theta^1_n \sin\bar \Theta^2_n \cos\bar \Theta^2_n (h+c_1\zeta_n\sqrt h), \\
\Theta^2_{n+1}&=\Theta^2_n+ R \sin\bar \Theta^1_n \left(
\Big{(}\frac{1}{I_1}-\frac{1}{I_3}\Big{)}
\cos^2 \bar \Theta^2_n
- \Big{(}\frac{1}{I_3}-\frac{1}{I_2}\Big{)}
\sin^2 \bar \Theta^2_n
\right)
(h+c_1\zeta_n\sqrt h),
\end{split}
\end{equation}
where $\bar \Theta^1_n=\frac{1}{2}(\Theta^1_n+\Theta^1_{n+1})$ and $\bar \Theta^2_n=\frac{1}{2}(\Theta^2_n+\Theta^2_{n+1})$.
Then by the mapping \eqref{sps-eq-9}, we obtain the following scheme for the original stochastic rigid body system \eqref{model;rigid}
\begin{align}
\label{snm;rg2}
Y^1_n = R \cos \Theta^1_n \cos \Theta^2_n,
\quad
Y^2_n = R \cos \Theta^1_n \sin \Theta^2_n, \quad
Y^3_n = R\sin \Theta^1_n,
\end{align}
which naturally satisfies 
\begin{align*}
(Y^1_n)^2+(Y^2_n)^2+(Y^3_n)^2=2\mathcal C_1,
\end{align*}
meaning that the scheme \eqref{snm;rg2} preserves the Casimir function.

Alternatively, we can also convert the SDE \eqref{sps-eq-10} to its equivalent It\^o form, and then use the Euler-Maruyama method or Milstein method, together with the spherical coordinate transformation, to construct numerical schemes preserving the Casimir function of the stochastic rigid body system \eqref{model;rigid}. In the following, we call numerical schemes resulted from the spherical transformation for the stochastic rigid body system \eqref{model;rigid} the `spherical schemes'. 

Next we derive the root mean-square convergence order of the spherical schemes. Denote a numerical scheme applied to \eqref{sps-eq-10}
by $\{\mathcal P_n,\mathcal Q_n\}_{n=0}^N,$ and its spherically transformed scheme for the original stochastic rigid body system \eqref{model;rigid} by $\{X_n^1,X_n^2,X_n^3\}_{n=0}^N$.

\begin{tm}\label{tmm}
If the scheme $\{\mathcal P_n,\mathcal Q_n\}_{n=0}^N$ applied to \eqref{sps-eq-10} is of root mean-square convergence order $k,$ then $\{X_n^1,X_n^2,X_n^3\}_{n=0}^N$ for \eqref{model;rigid} is also of root mean-square convergence order $k,$  that is, for any $T>0$ with $0=t_0<t_1<\cdots<t_N=T$, $h=t_{j+1}-t_j$ ($j=0,\cdots, N-1$),
\begin{equation}
\label{tm;eq}
\E\left\|\left(y_1(T),y_2(T),y_3(T)\right)-(X_N^1,X_N^2,X_N^3)\right\|^2=O(h^{2k}).
\end{equation}
\end{tm}
{\bf Proof.}
Since the coefficients of \eqref{sps-eq-10} are globally Lipschitz continuous, the scheme $\{\mathcal P_n,\mathcal Q_n\}_{n=0}^N$ which is of root mean-square convergence order $k$ has finite moments. Using the coordinate transformation $\phi,$ the left hand side of $\eqref{tm;eq}$ can be written as
\begin{align*}
\mathcal I=\E\|\phi(\theta_1(T),\theta_2(T))-\phi(\mathcal P_N,\mathcal Q_N)\|^2.
\end{align*}
Using the Lipschitz continuity of the mapping $\phi$, we have
\begin{align*}
\mathcal I\le&\E\|\phi(\theta_1(T),\theta_2(T))-\phi(\theta_1(T),\mathcal Q_N)\|^2
+\E\|\phi(\theta_1(T),\mathcal Q_N))-\phi(\mathcal P_N,\mathcal Q_N)\|^2\\
\leq& K\E\|\theta_2(T)-\mathcal Q_N\|^2
+K\E\|\theta_1(T)-\mathcal P_N\|^2
=KO(h^{2k}),
\end{align*}
where $K$ is a sufficiently large number independent of $h.$
\hfill$\square$\\

It follows from Theorem \ref{tmm} that the numerical scheme \eqref{snm;rg2} has root mean-square convergence order 1, since the midpoint rule \eqref{snm;sym2} is of root mean-square order 1.
\begin{rk} \hfill
\begin{itemize}
\item We can see from the proof above that, if the system resulted from a coordinate transformation $\varphi$, 
e.g. the system \eqref{sps-eq-10}, has globally Lipschitz continuous coefficients, and the transformation $\varphi$ is Lipschitz continuous, then the numerical schemes before and after the inverse transformation $\varphi^{-1}$, for the transfromed system and the original SPS, respectively, have the same root mean-square convergence order. 
\item For the canonical transformation \eqref{trans1}, however, we see that it is not globally Lipschitz continuous, which causes difficulties for theoretical  analysis on the root mean-square convergence order of the numerical scheme \eqref{snm;rg1} arising from this coordinate transformation. This is also the case for the scheme \eqref{lvscheme} in Section \ref{sec4.2}. We will then illustrate empirical analysis of the mean-square order of \eqref{snm;rg1} and \eqref{lvscheme} via numerical tests.
\end{itemize}
\end{rk}
\subsubsection{The three-dimensional stochastic Lotka-Volterra system}\label{sec4.2}
Consider the stochastic Lotka-Volterra system 
\begin{align}
\label{model:3DLV}
dy=B_2(y)\nabla K_2(y)(dt+c_2\circ dW(t)), \quad y(0)=(y_1^0,y_2^0,y_3^0)^\top,
\end{align}
where $y_i^0>0$ $(i=1,2,3)$, 
$K_2=aby_1+y_2-ay_3+v\ln y_2-\mu \ln y_3$,
$$B_2=
\begin{bmatrix}
0& ry_1y_2& bry_1y_3\\
-ry_1y_2& 0& y_2y_3\\
-bry_1y_3& -y_2y_3& 0
\end{bmatrix},
$$ and $a,b,c_2,r,v,\mu$ are constants.
As described in Section \ref{lvexa}, its solution is positive for all $t$ almost surely, and the Casimir function is 
\begin{equation}\label{casi2}
C_2 (y) = \frac{1}{r}\ln y_1 - b \ln y_2 + \ln y_3 \equiv \frac{1}{r}\ln y_1^0 - b \ln y_2^0 + \ln y_3^0=:\mathcal C_2.
\end{equation}

Analogous to the procedure for the stochastic rigid body system, we first look for a canonical coordinate transformation 
  $\bar{y} = (\bar{y}_1(y), \bar{y}_2(y), \bar{y}_3(y))^\top$ with Jacobian matrix $A(y)=\frac{\partial \bar y}{\partial y}$ satisfying
\begin{align}\label{canna}
A(y) B_2 (y) A(y)^{\mathrm{T}} = B_0:=
\begin{bmatrix}
0 & 1 & 0 \\ -1 & 0 & 0 \\ 0 & 0 & 0
\end{bmatrix}.
\end{align}
Solving the partial differential equation systems \eqref{canna} based on the Poisson bracket defined by $B_2$, we find the following coordinate transformation
\begin{equation*}
\bar{y}_1 = \ln y_3, \quad \bar{y}_2 = - \ln y_2, \quad \bar{y}_3 = \mathcal C_2,
\end{equation*}
and its inverse 
\begin{equation}
\label{map2}
y_1 = \exp (r \left( \mathcal C_2-\bar{y}_1 - b \bar{y}_2 \right) ), \quad y_2 = \exp (-\bar{y}_2), \quad y_3 = \exp (\bar{y}_1).
\end{equation}
Then, denoting $(\bar y_1,\bar y_2)=(P,Q)$, we get the following SHS driven by multiplicative noise
\begin{align}
\label{model:ESHS}
d\begin{bmatrix}
P\\
Q
\end{bmatrix}=
\begin{bmatrix}
0& -1\\
1& 0
\end{bmatrix}
\nabla H(P,Q)(dt+c_2\circ dW(t)),
\end{align}
where $H(P,Q)=-ab\exp(r(\mathcal C_2-P-bQ))
-\exp(-Q)+vQ+a\exp(P)+\mu P.$
More explicitly, 
\begin{align*}
d\begin{bmatrix}
P\\
Q
\end{bmatrix}=
\begin{bmatrix}
0& -1\\
1& 0
\end{bmatrix}
\begin{bmatrix}
abr\exp (r(\mathcal C_2-P-bQ))
+a\exp(P)+\mu\\
ab^2r\exp (r(\mathcal C_2-P-bQ))
+\exp(-Q)+v
\end{bmatrix}
(dt+c_2\circ dW(t)).
\end{align*}
Applying the symplectic scheme \eqref{sns1} to \eqref{model:ESHS}, we get  
\begin{align}\label{lvscheme}
P_{n+1}&=P_n-(h+c_2\zeta_n\sqrt{h})\left[ab^2r\exp(r(\mathcal C_2-\bar P_n-b\bar Q_n))+\exp(-\bar Q_n)+v\right]\nonumber\\
&-c_{\alpha}\zeta_n^2h\left[-2a^2b^4r^3\exp(2r(\mathcal C_2-\bar P_n-b\bar Q_n))\right.\nonumber\\
&-(rb+1)abr\exp(r(\mathcal C_2-\bar P_n-b\bar Q_n)-\bar Q_n)\nonumber\\
&-rb(vabr+\mu ab^2r)\exp(r(\mathcal C_2-\bar P_n-b\bar Q_n))\nonumber\\
&-a^2b^3r^2\exp(r(\mathcal C_2-\bar P_n-b\bar Q_n)+\bar P_n)\nonumber\\
&\left.-a\exp(\bar P_n-\bar Q_n)-\mu \exp(-\bar Q_n)\right],\nonumber\\
Q_{n+1}&=Q_n+(h+c_2\zeta_n\sqrt{h})\left[abr\exp(r(\mathcal C_2-\bar P_n-b\bar Q_n))+a\exp(\bar P_n)+\mu\right]\nonumber\\
&+c_{\alpha}\zeta_n^2h\left[-2a^2b^3r^3\exp(2r(\mathcal C_2-\bar P_n-b\bar Q_n))\right.\nonumber\\
&-abr^2\exp(r(\mathcal C_2-\bar P_n-b\bar Q_n)-\bar Q_n)\nonumber\\
&-(vabr^2+\mu ab^2r^2)\exp(r(\mathcal C_2-\bar P_n-b\bar Q_n))\nonumber\\
&+(1-r)a^2b^2r\exp(r(\mathcal C_2-\bar P_n-b\bar Q_n)+\bar P_n)\nonumber\\
&\left.+a\exp(\bar P_n-\bar Q_n)+av\exp (\bar P_n)\right],
\end{align}
where $c_{\alpha}=c_2^2(\alpha-\frac{1}{2})$, $\bar P_n=(1-\alpha)P_n+\alpha P_{n+1}$, $\bar Q_n=\alpha Q_n+(1-\alpha)Q_{n+1}$, $\alpha\in [0,1]$. Then, by using the inverse transformation \eqref{map2}, we obtain the following numerical scheme for the original stochastic Lotka-Volterra system \eqref{model:3DLV}
\begin{equation}
\label{lvscheme}
Y^1_{n} = \exp (r \left( \mathcal C_2-P_{n} - b Q_{n} \right) ), \quad Y^2 _{n}= \exp (-Q_{n}), \quad Y^3_{n} = \exp (P_{n}).
\end{equation}

Obviously, the $\alpha$-generating schemes \eqref{lvscheme} preserve the positivity of the solution of the original stochastic Lotka-Volterra system. Moreover, it can be easily verified that they also preserve the Casimir function \eqref{casi2} of the system. Using the proof for Theorem \ref{gen1}, we can show that the $\alpha$-generating schemes \eqref{lvscheme} preserve the Poisson structure of the stochastic Lotka-Volterra system \eqref{model:3DLV}. 

\section{Numerical illustrations}\label{nmil}
\subsection{The stochastic rigid body system}
In this subsection we demonstrate the numerical behavior of the $\alpha$-generating schemes \eqref{snm;rg1}, and that of the spherical scheme \eqref{snm;rg2} for the stochastic rigid body system \eqref{model;rigid}.
\begin{figure}[htbp]
\centering
\subfigure[$\alpha=0$]{\label{f1a}
\includegraphics[width=5.5cm,height=5cm]{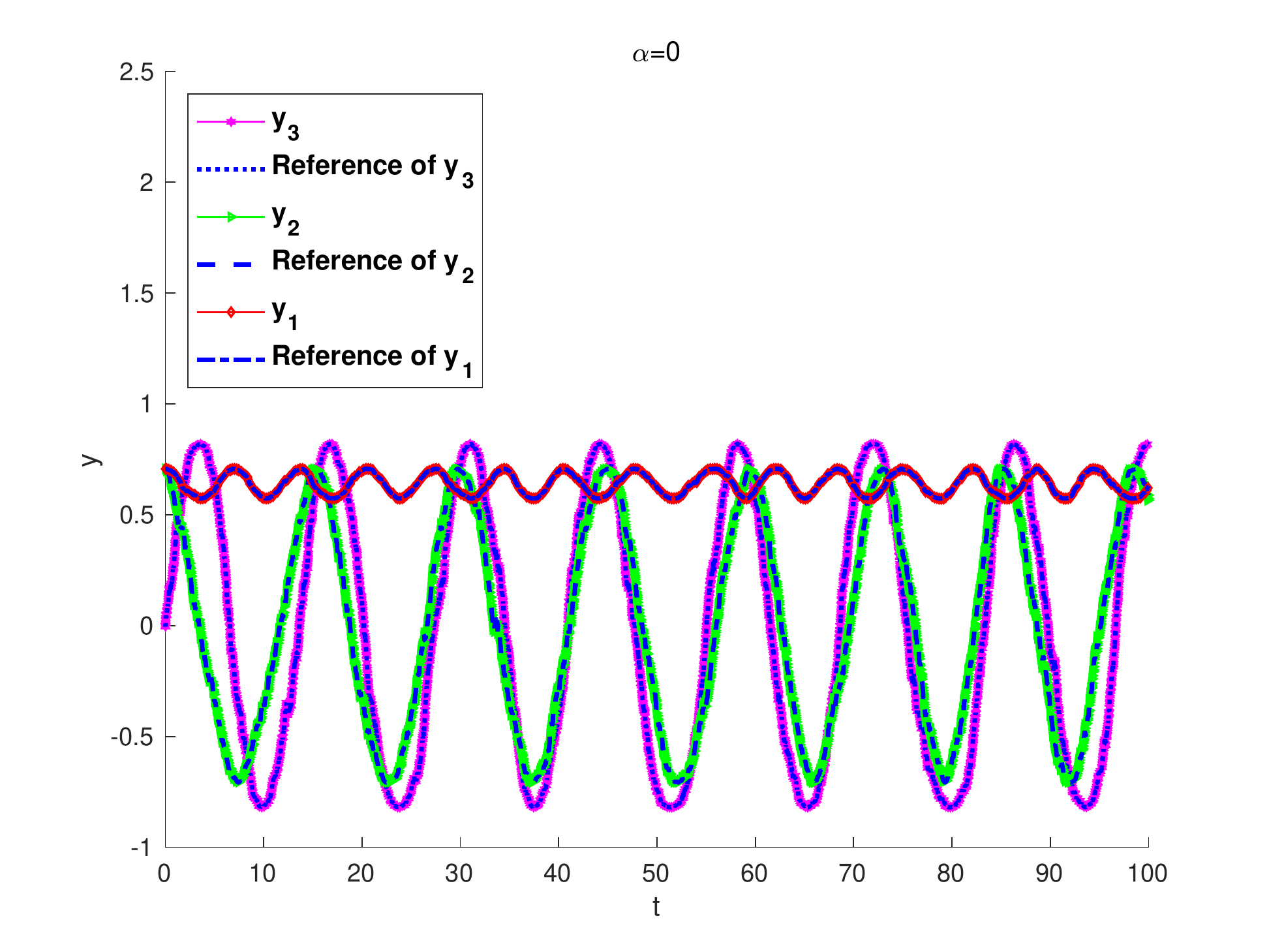}}
\subfigure[$\alpha=1$]{\label{f1b}
\includegraphics[width=5.5cm,height=5cm]{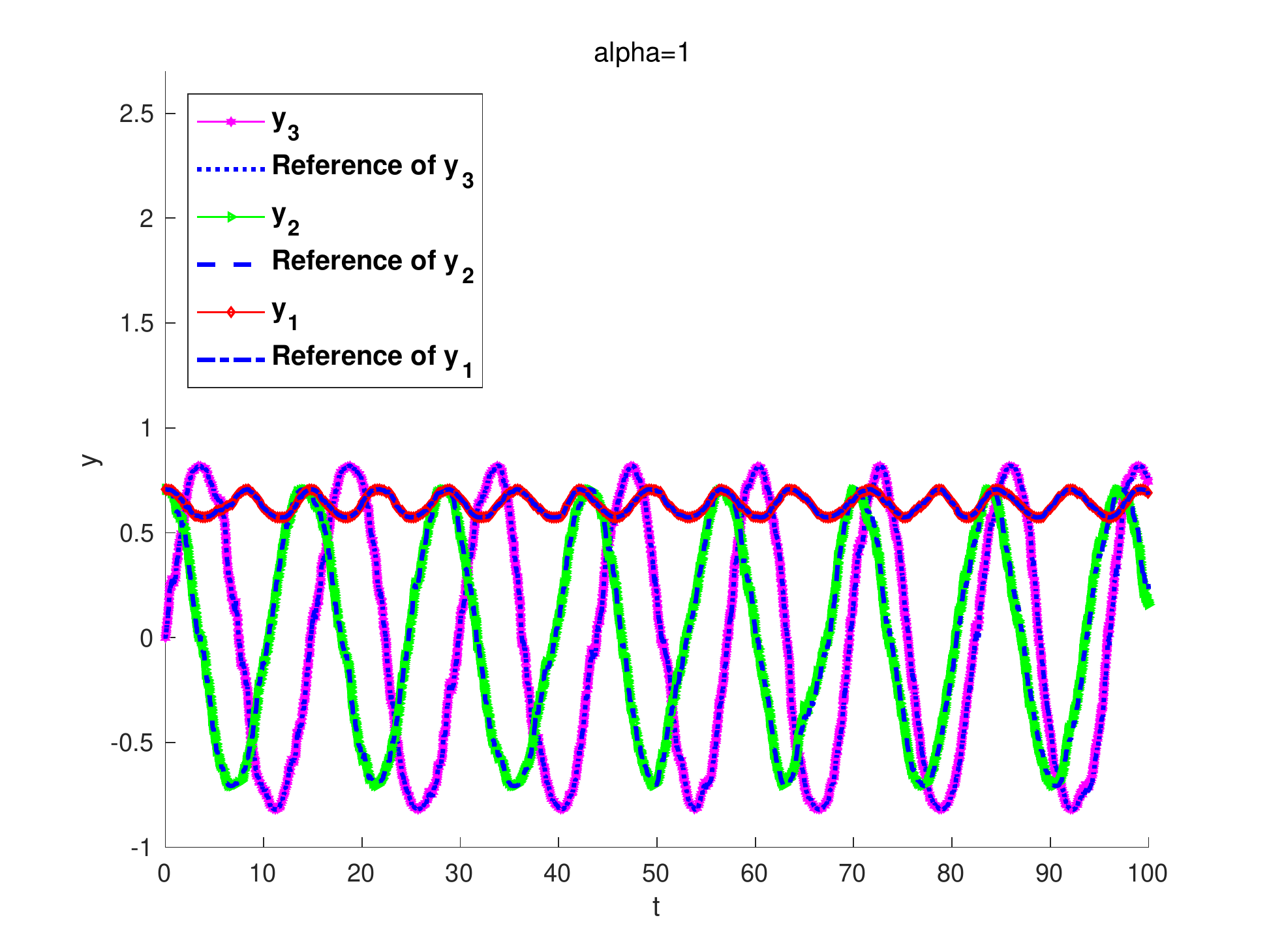}}
\subfigure[$\alpha=0.5$]{\label{f1c}
\includegraphics[width=5.5cm,height=5cm]{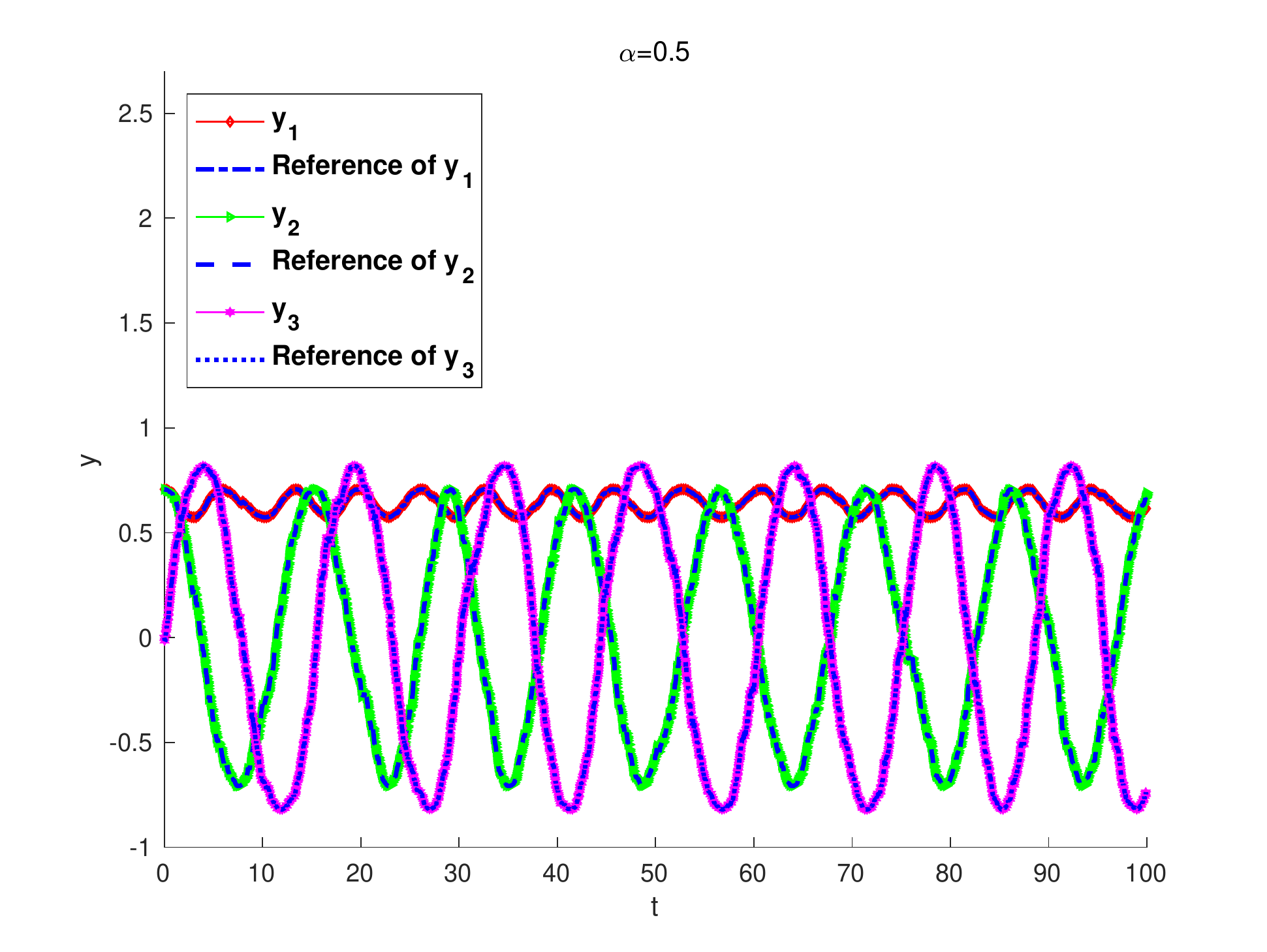}}
\subfigure[spherical scheme]{\label{f1d}
\includegraphics[width=5.5cm,height=5cm]{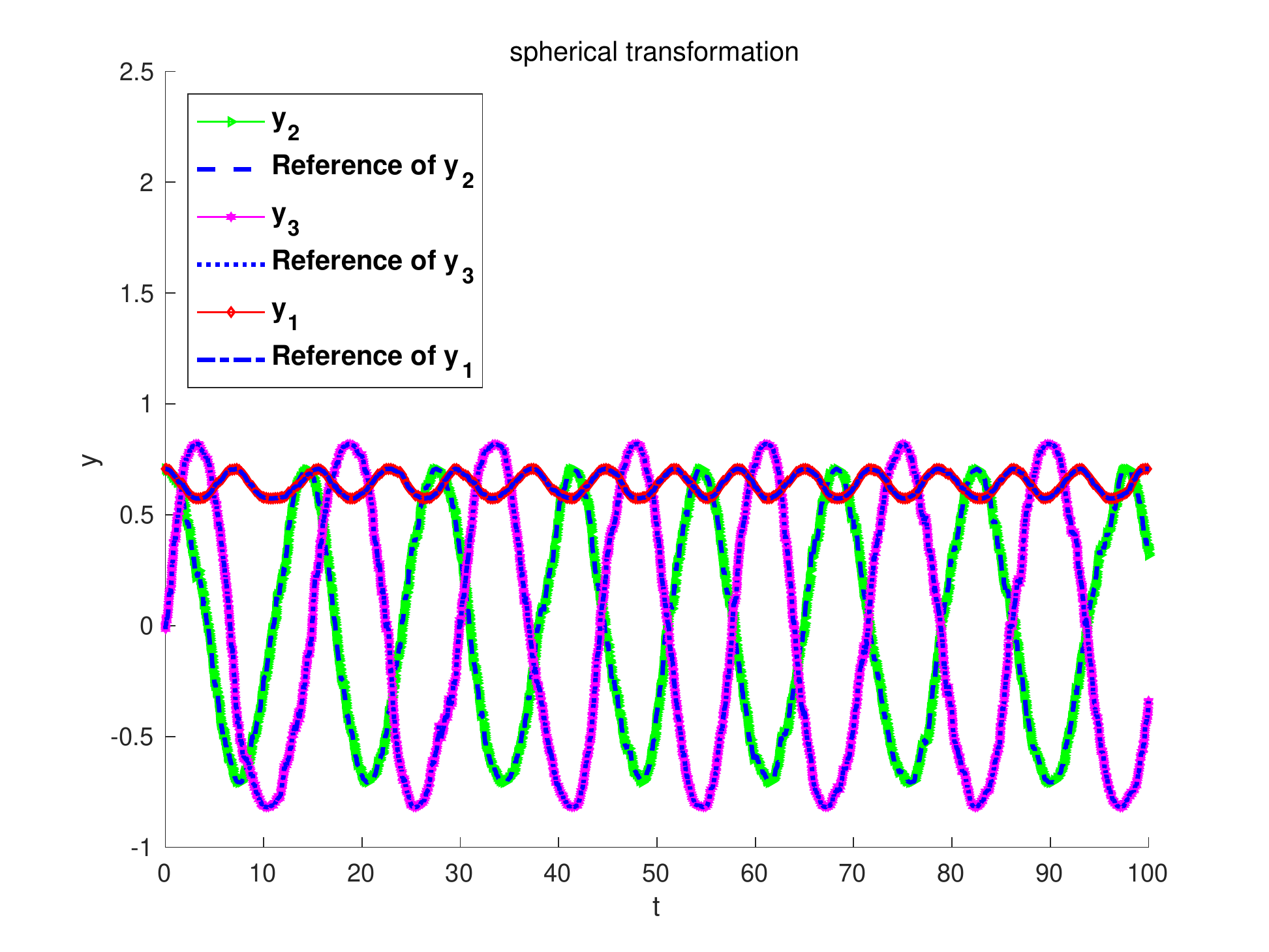}}
\caption{Sample paths of $y_1$, $y_2$ and $y_3$ produced by the $\alpha$-generating schemes and the spherical scheme }
\label{f1}
\end{figure}

Figure \ref{f1} shows the sample paths of $y_1$, $y_2$ and $y_3$ of the stochastic rigid body system \eqref{model;rigid} produced by the $\alpha$-generating schemes \eqref{snm;rg1} with $\alpha=0$ (Figure \ref{f1a}), $\alpha=1$ (Figure \ref{f1b}), $\alpha=0.5$ (Figure \ref{f1c}), and by the spherical scheme \eqref{snm;rg2} (Figure \ref{f1d}). The reference solutions of $y_1$, $y_2$ and $y_3$ (blue) are approximated by midpoint rule with time step $10^{-5}$. The constants take the value $I_1=\sqrt{2}+\sqrt{\frac{2}{1.51}}$, $I_2=\sqrt{2}-0.51\sqrt{\frac{2}{1.51}}$, $I_3=1$, $c_1=0.2$. The initial values of $y$ are $y_1^0=y_2^0=\frac{1}{\sqrt{2}}$, $y_3^0=0$. We take time step $h=0.01$, $err=10^{-12}$ as the error bound for stopping the inner iterations within each time step by implementing the implicit schemes. We can see that, all the numerical sample paths coincide very well with the reference solutions.

\begin{figure}[htbp]
\centering
\subfigure[Casimir by $\alpha=0,0.5,1$, and the spherical scheme]{\label{f2a}
\includegraphics[width=5.5cm,height=5cm]{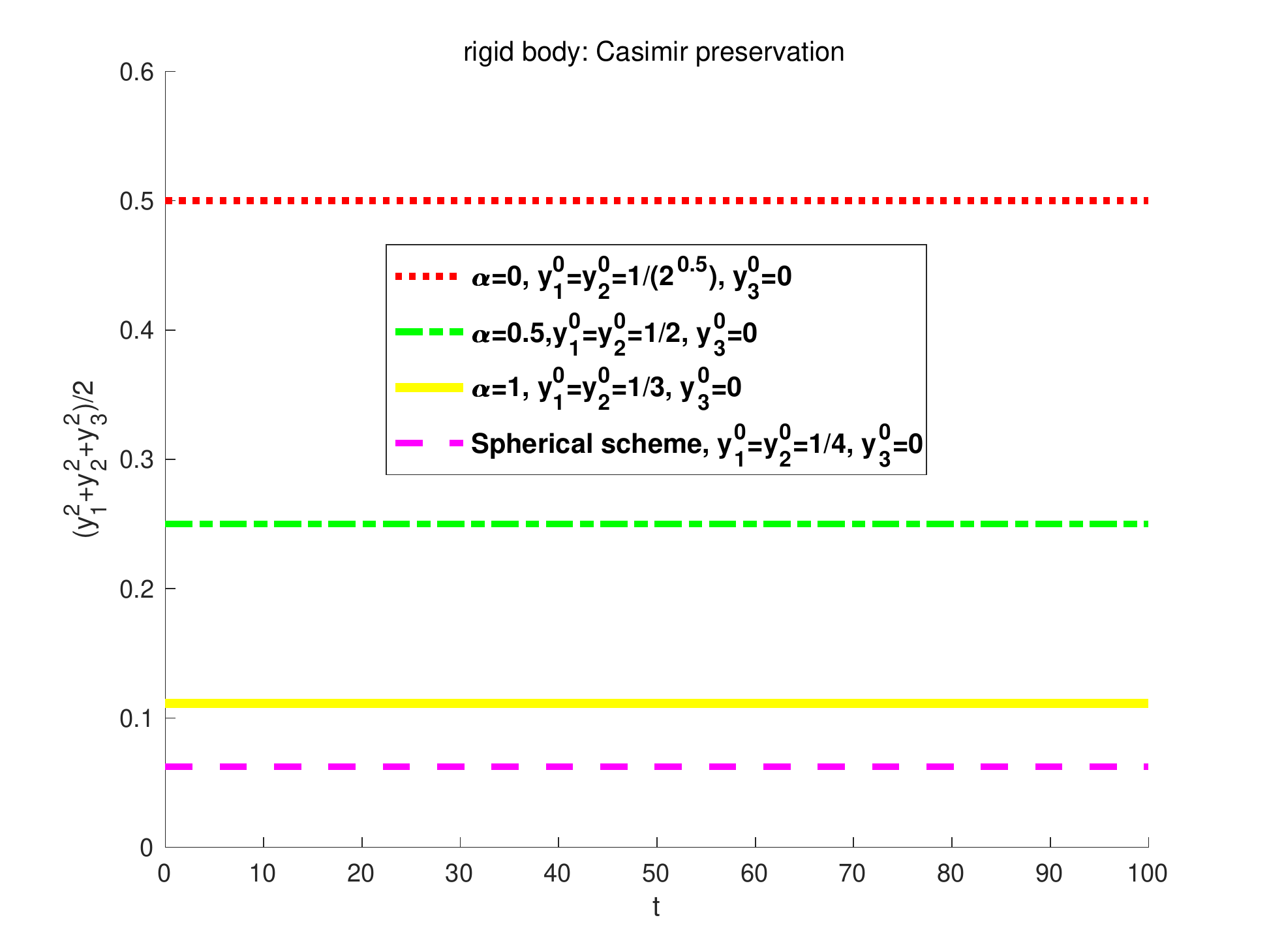}}
\subfigure[Comparison with the Euler-Maruyama scheme ]{\label{f2b}
\includegraphics[width=5.5cm,height=5cm]{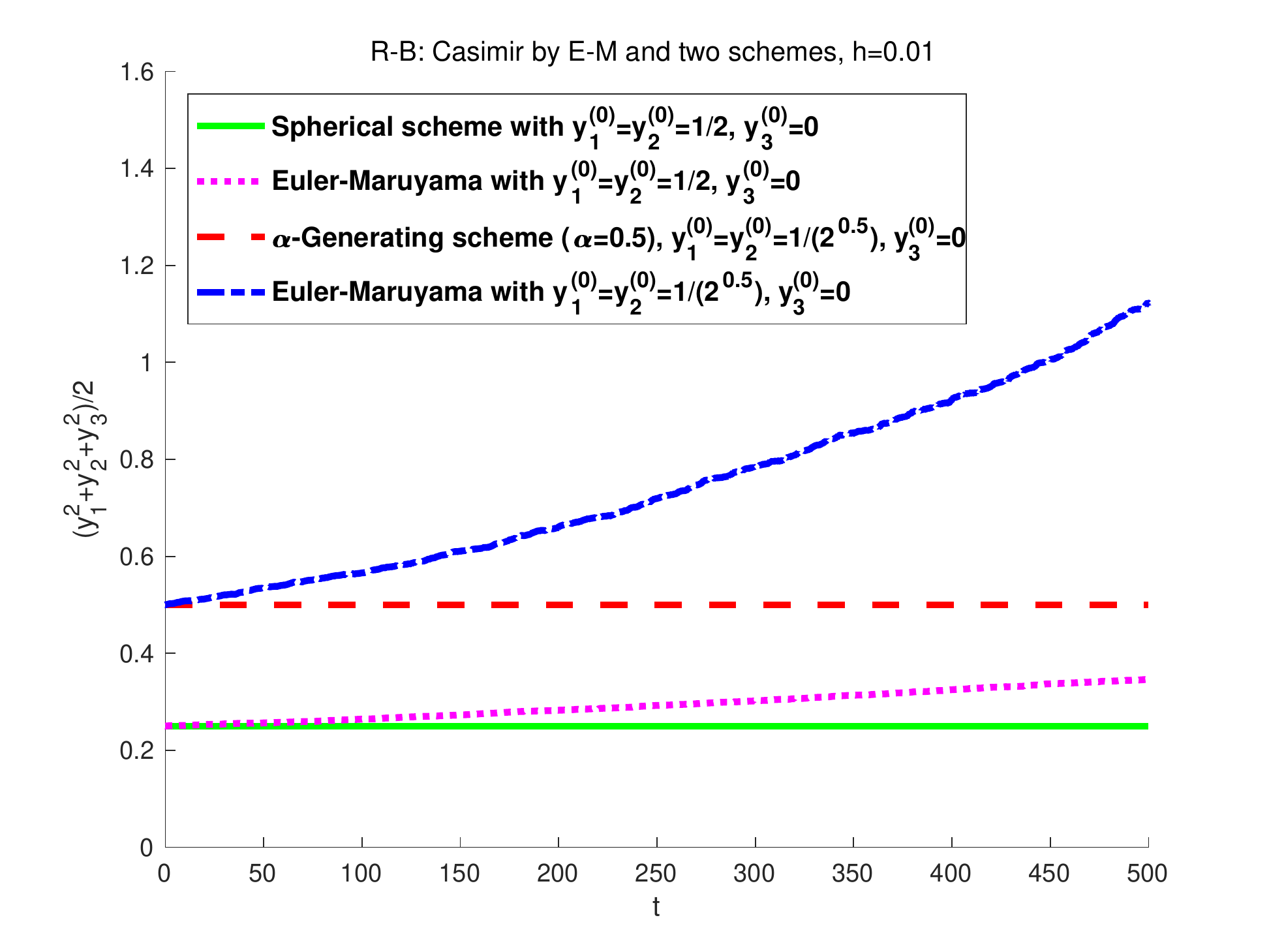}}
\caption{Numerical evolution of the Casimir function by the $\alpha$-generating schemes, the spherical scheme, and the Euler-Maruyama scheme }
\label{f2}
\end{figure}

Figure \ref{f2} illustrates the evolution of the Casimir produced by the $\alpha$-generating schemes with $\alpha=0,0.5,1$, and the spherical scheme, each with a different initial value of $(y_1^0,y_2^0,y_3^0)$ giving different Casimir values $\mathcal C_1$. Figure \ref{f2a} shows clearly the exact preservation of the Casimir function by the $\alpha$-generating schemes and the spherical scheme. Figure \ref{f2b} compares the Casimir evolution by the $\alpha$-generating schemes and the spherical scheme with that by the Euler-Maruyama scheme. We can see that the Euler-Maruyama scheme fails to preserve the Casimir function. Note that by using the Euler-Maruyama scheme, for consistency, we need to transform the system \eqref{model;rigid} to its equivalent It\^o SDE. In Figure \ref{f2a}, for $\alpha=0$, we take $y_1^0=y_2^0=\frac{1}{\sqrt{2}}$, $y_3^0=0$, that is $\mathcal C_1=\frac{1}{2}$. For $\alpha=0.5,$ $y_1^0=y_2^0=\frac{1}{2}$, $y_3^0=0$, then $\mathcal C_1=\frac{1}{4}$. For $\alpha=1$, $y_1^0=y_2^0=\frac{1}{3}$, $y_3^0=0$, then $\mathcal C_1=\frac{1}{9}$. For the spherical scheme, $y_1^0=y_2^0=\frac{1}{4}$, $y_3^0=0$, then $\mathcal C_1=\frac{1}{16}$. T=100. In Figure \ref{f2b}, for comparison between the $\alpha=0.5$-generating scheme and the Euler-Maruyama scheme, we let $y_1^0=y_2^0=\frac{1}{\sqrt{2}}$, $y_3^0=0$, and for that between the spherical and the Euler-Maruyama scheme, we let $y_1^0=y_2^0=\frac{1}{2}$, $y_3^0=0$. $T=500$. In both subfigures, $h=0.01$. Other data are the same with those for Figure \ref{f1}.

\begin{figure}[htbp]
\centering
\subfigure [$\alpha=0,0.5,1$]{\label{f3a}
\includegraphics[width=5.5cm,height=5cm]{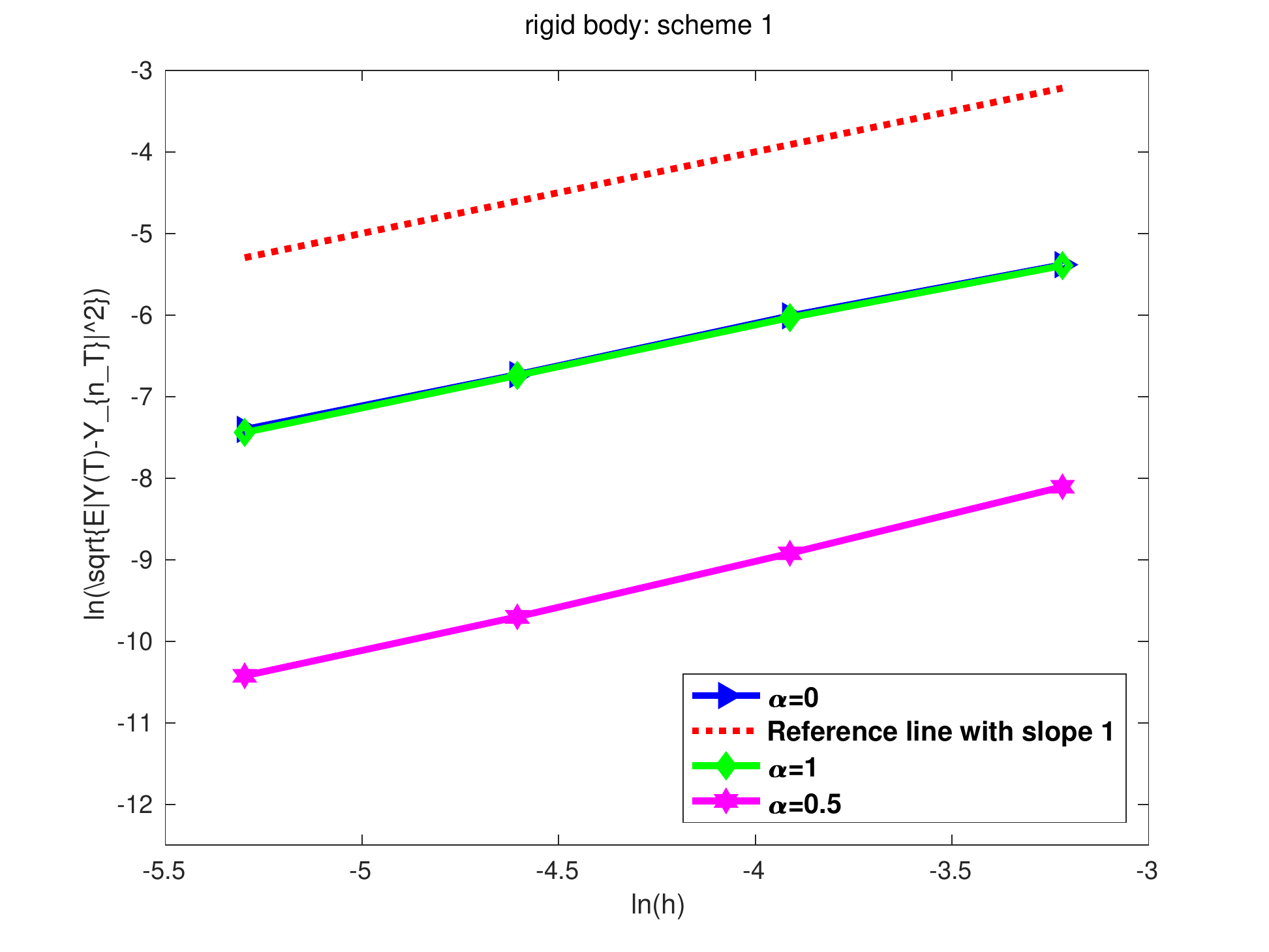}}
\subfigure[The spherical scheme ]{\label{f3b}
\includegraphics[width=5.5cm,height=5cm]{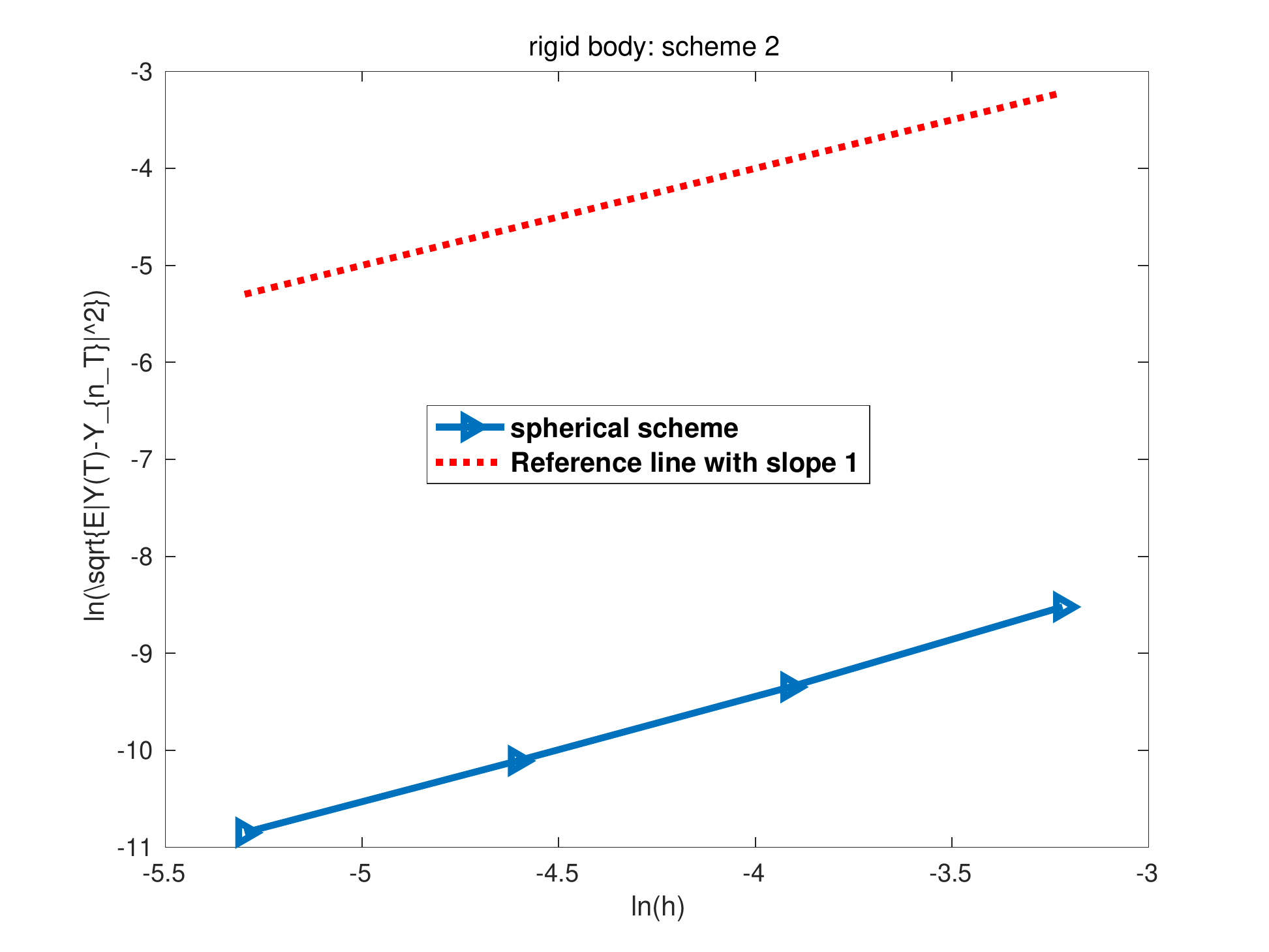}}
\caption{Mean-square order of the $\alpha$-generating schemes and the spherical scheme }
\label{f3}
\end{figure}

Figure \ref{f3} shows the root mean-square convergence order of the $\alpha$-generating schemes with $\alpha=0,0.5,1$, and that of the spherical scheme. From Figure \ref{f3a} we see that the $\alpha$-generating schemes have root mean-square order 1, and the $\alpha=0.5$-generating scheme is with smaller error than $\alpha=0$ and $0.5$, and the lines for $\alpha=0$ and $1$ coincide visually. We can see from Figure \ref{f3b} that the spherical scheme is also of root mean-square convergence order 1. In both subfigures we take $h=[0.005,0.01,0.02,0.04]$, $T=10$, and 500 samples for approximating the expectations. Other data are the same with those for Figure \ref{f1}.

\subsection{The stochastic Lotka-Volterra system}
In this section we observe the behavior of  $\alpha$-generating schemes \eqref{lvscheme} for the stochastic Lotka-Volterra system \eqref{model:3DLV} via numerical experiments.
\begin{figure}[htbp]
\centering
\subfigure[$\alpha=0$]{\label{f4a}
\includegraphics[width=3.6cm,height=5cm]{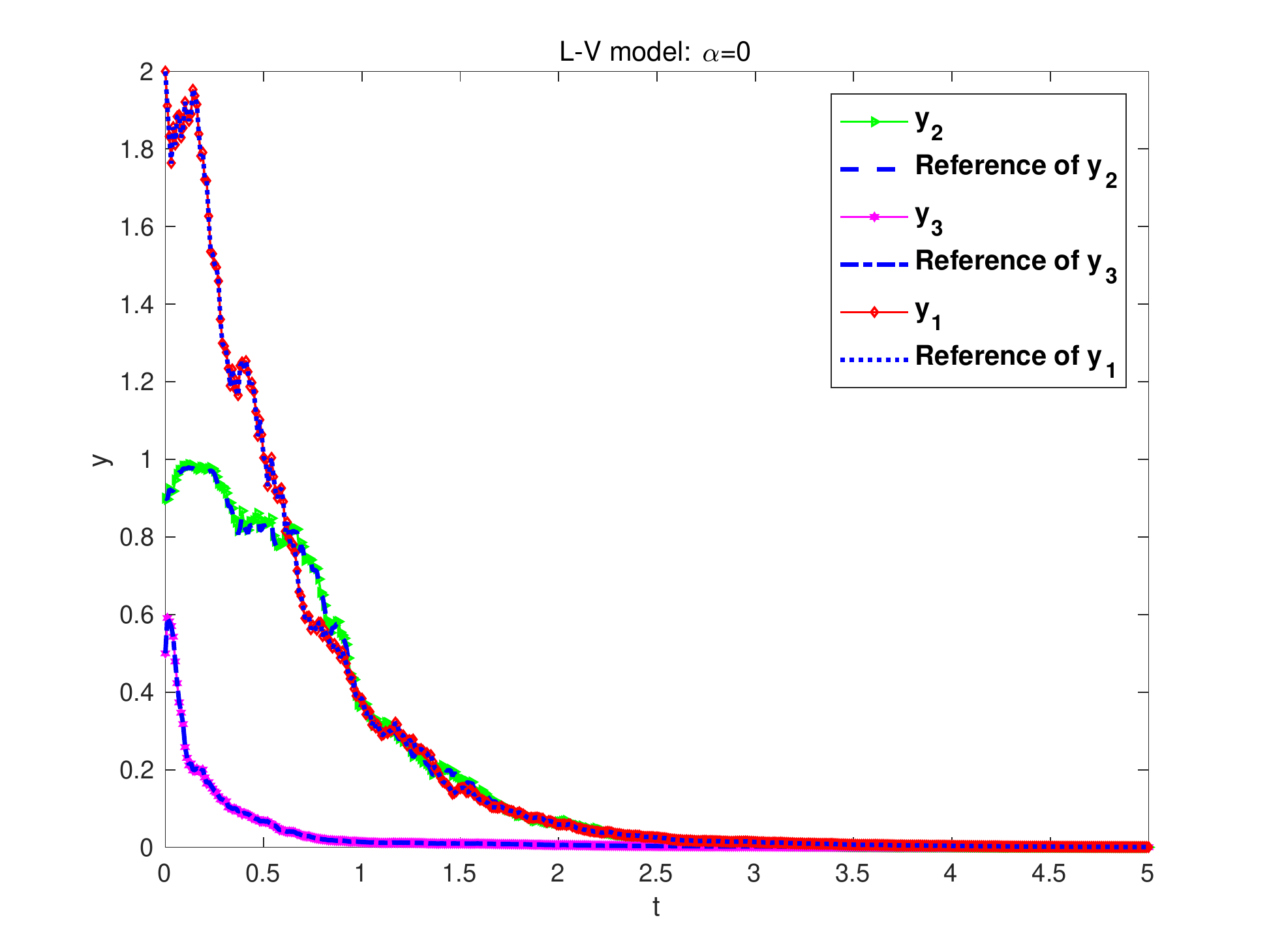}}
\subfigure[$\alpha=1$]{\label{f4b}
\includegraphics[width=3.6cm,height=5cm]{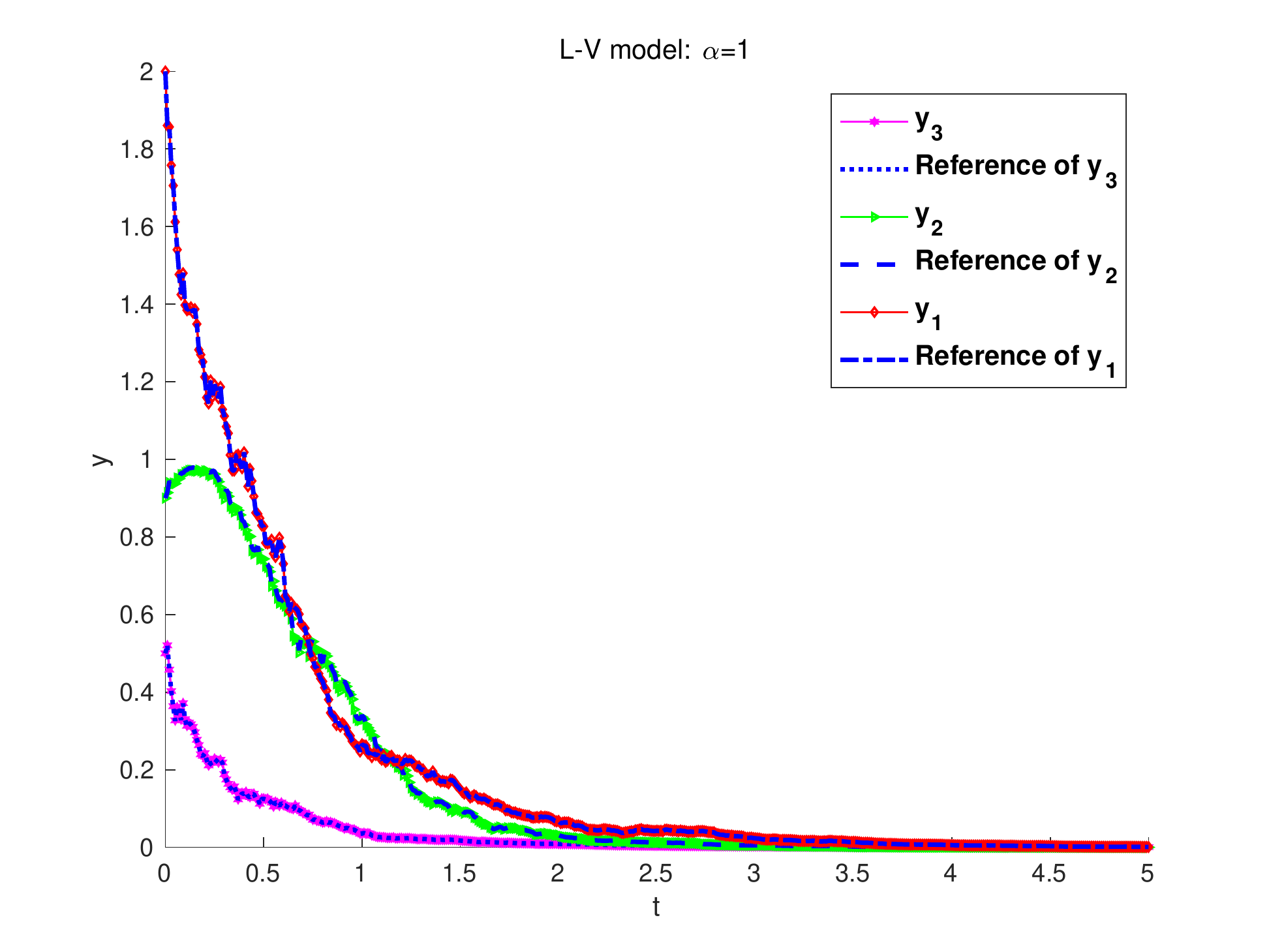}}
\subfigure[$\alpha=0.5$]{\label{f4c}
\includegraphics[width=3.6cm,height=5cm]{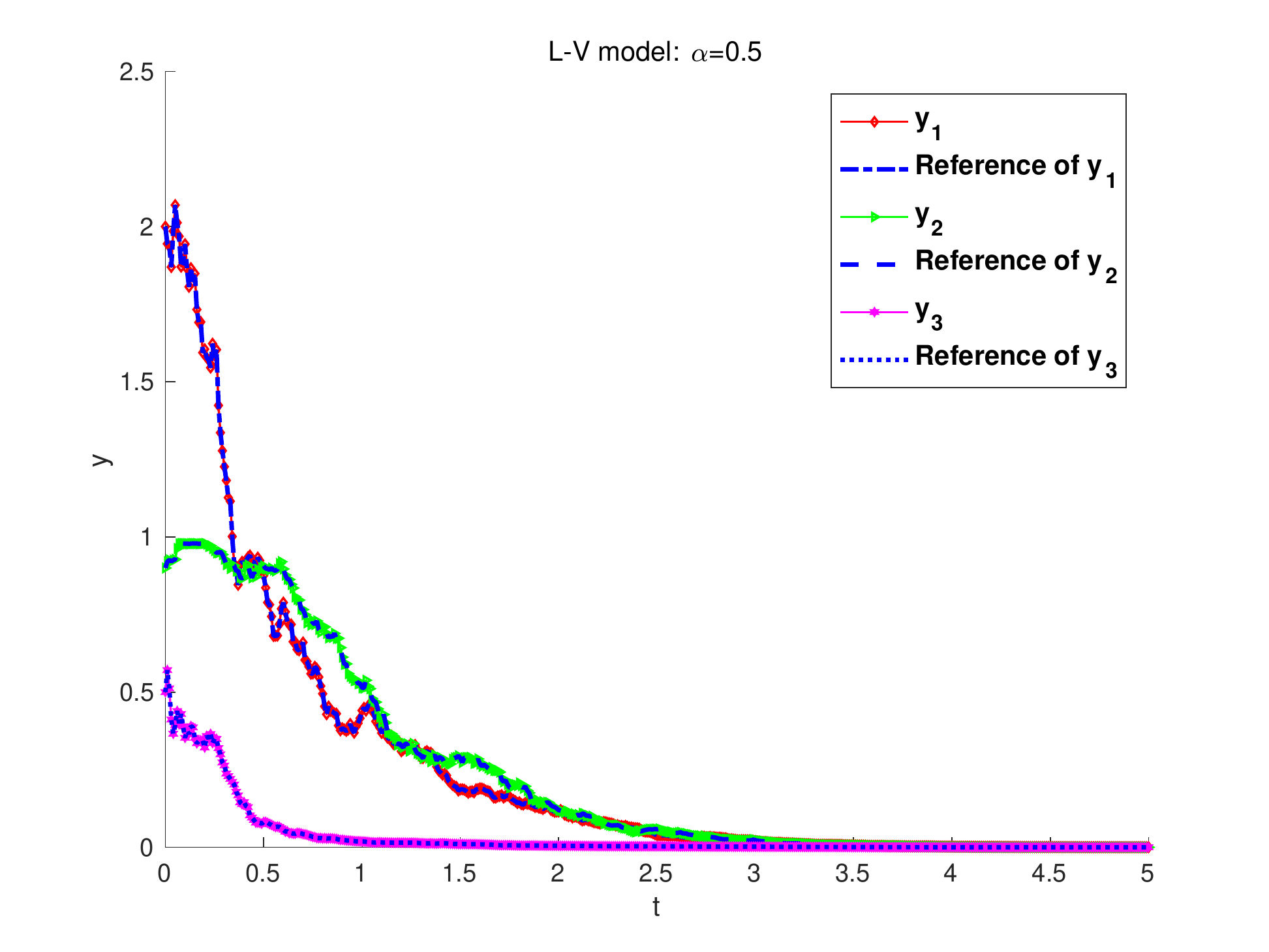}}
\caption{Sample paths of $y_1$, $y_2$ and $y_3$ produced by the $\alpha$-generating schemes \eqref{lvscheme} }
\label{f4}
\end{figure}

Figure \ref{f4} shows the sample paths of $y_1,y_2,y_3$ produced by the scheme \eqref{lvscheme} with $\alpha=0$ (Figure \ref{f4a}), $\alpha=1$ (Figure \ref{f4b}) and $\alpha=0.5$ (Figure \ref{f4c}). The reference solutions are simulated by midpoint rule with time step $10^{-5}$. The time step used in the numerical schemes is $h=0.01$. The initial data are $y_1^0=2$, $y_2^0=0.9$, $y_3^0=0.5$. The constants are $a=-2$, $b=-1$, $c_2=0.2$, $r=-0.5$, $\mu=2$, $v=1$. We can observe very good coincidence between the numerical and the reference solutions.

\begin{figure}[htbp]
\centering
\subfigure [Comparison of Casimir evolution]{\label{f5a}
\includegraphics[width=5.5cm,height=5cm]{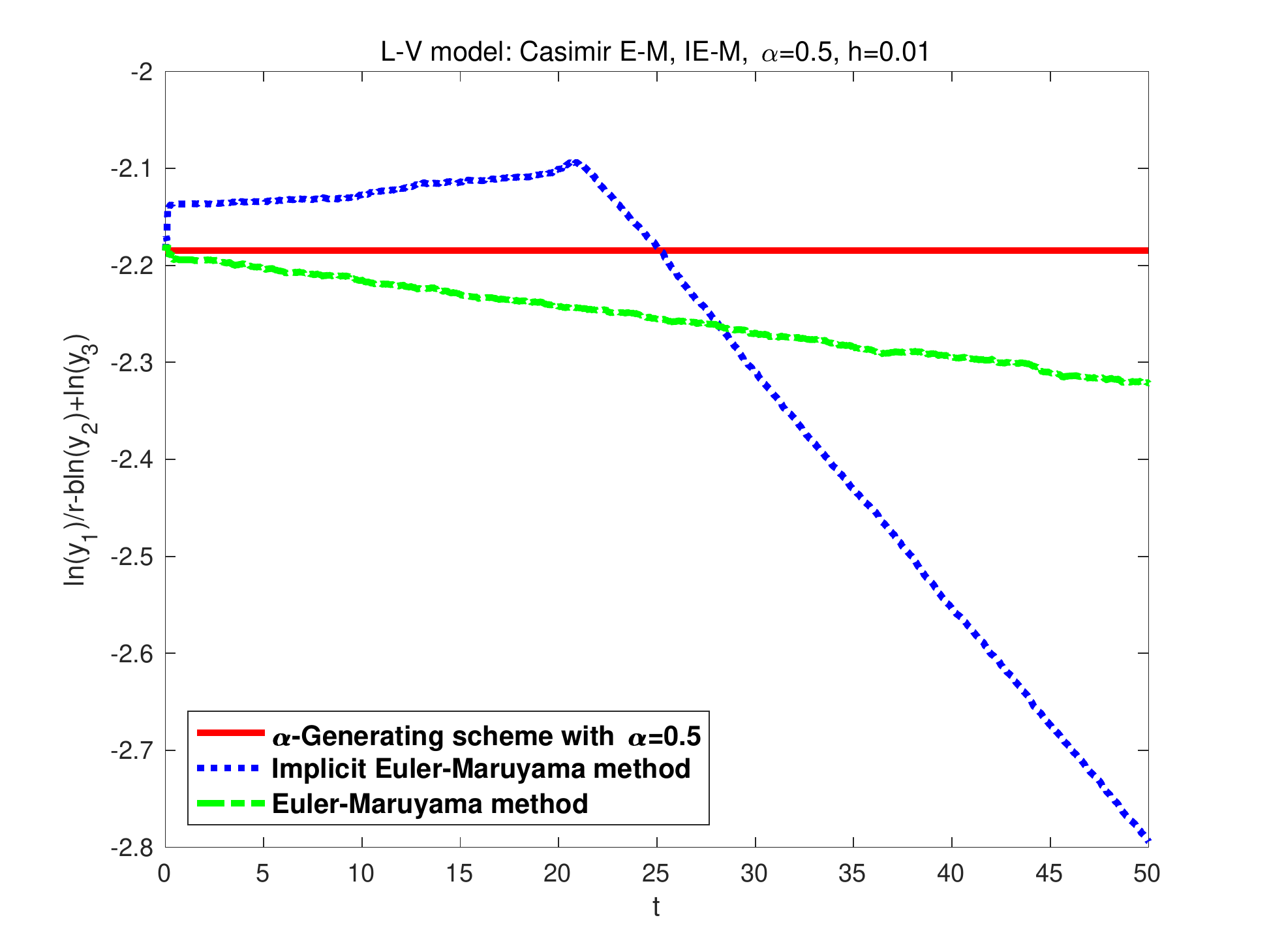}}
\subfigure[Mean-square orders  ]{\label{f5b}
\includegraphics[width=5.5cm,height=5cm]{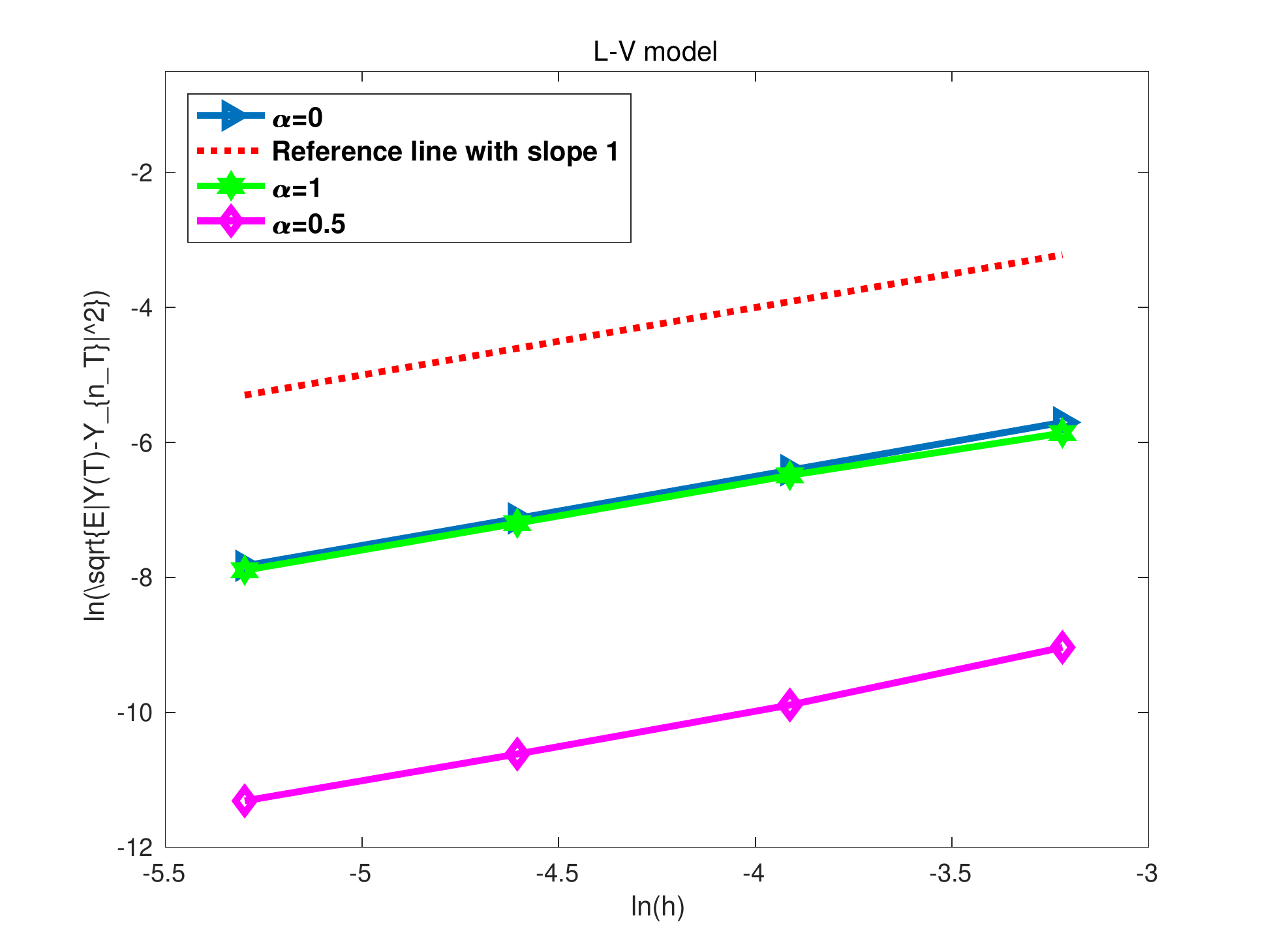}}
\caption{Casimir evolution and mean-square order of the $\alpha$-generating schemes \eqref{lvscheme} }
\label{f5}
\end{figure}

Figure \ref{f5a} compares the numerical evolution of the Casimir function produced by the scheme \eqref{lvscheme} with $\alpha=0.5$ (red solid) and by the explicit (green dash-dotted) and implicit Euler-Maruyama (blue dotted) schemes. Obviously the explicit and implicit Euler-Maruyama methods can not preserve the Casimir function, while the $\alpha$-generating scheme \eqref{lvscheme} can. The data are the same with those for Figure \ref{f4}. 

Figure \ref{f5b} shows the root mean-square convergence order of the scheme \eqref{lvscheme} with $\alpha=0,1$ and $0.5$. We see that they are of root mean-square convergence order 1, and the $\alpha=0.5$ scheme (pink) has the smallest error among the three. We take $h=[0.005,0.01,0.02,0.04]$ for plotting the lines, $T=2$, and 500 samples for approximating the expectation. Other data are the same with those for Figure \ref{f4}.



\section{Conclusions}
The proposed numerical algorithms based on the Darboux-Lie theorem and the $\alpha$-generating function approach for stochastic Poisson systems of arbitrary dimensions are proved to be efficient and structure-preserving for both the Poisson structure and the Casimir functions. It provides a large variety of stochastic Poisson integrators. 

Beyond the efficiency of the numerical methods, verified theoretically and numerically, we also emphasize the flexibility of the algorithm, in that it allows free choice of $\alpha\in[0,1]$,  and different canonical coordinate transformations. 

\section*{Acknowledgments}
Authors are funded by the National Natural Science Foundation of China (No. 91530118, No. 91130003, No. 11021101, No. 11290142, No. 11971458, No. 11471310, No.11071251).

\bibliographystyle{plain}
\bibliography{references}


\end{document}